\input ./style/arxiv-general.cfg
\documentclass[MSNbibl,number,citesort,seceqn,dvips]{arxbj}
\makeatletter
   \@ifpackageloaded{graphicx}{}{\usepackage{graphicx}}
\makeatother
\usepackage{mathrsfs}
\usepackage{mathbh}

\aid{0}
\volume{21}
\issue{4}
\pubyear{2015}
\firstpage{2430}
\lastpage{2456}
\doi{10.3150/14-BEJ650} 
\docsubty{FLA}

\makeatletter

\newcommand{\rrvert}{\vert}
\newcommand{\rrVert}{\Vert}
\newcommand{\llvert}{\vert}
\newcommand{\llVert}{\Vert}
\newcommand{\eqref}[1]{(\ref{#1})}
\newtheorem{lemma}{Lemma}
\newtheorem{theorem}{Theorem}
\newremark{Example}{Example}

\newcommand{\goth}{\mathfrak}
\newcommand{\Pb}{\mathbf{P}}
\newcommand{\Ex}{\mathbf{E}}
\renewcommand{\SS}{\mathbb{S}}
\newcommand{\ZZ}{\mathbb{Z}}
\newcommand{\RR}{\mathbb{R}}
\newcommand{\NN}{\mathbb{N}}
\newcommand{\II}{\mathbb{I}}
\newcommand{\JJ}{\mathbb{J}}
\newcommand{\one}{\mathbh{1}}
\makeatother

\begin{document}
\begin{frontmatter}

\title{On ADF goodness-of-fit tests for perturbed dynamical systems}
\runtitle{ADF GoF tests for dynamical systems}

\begin{aug}
\author[A]{\inits{Y.A.}\fnms{Yury A.}~\snm{Kutoyants}\corref{}\ead[label=e1]{kutoyants@univ-lemans.fr}}
\address[A]{Laboratoire de Statistique et Processus, Universit\'e du
Maine, 72085
Le Mans, France and Laboratory of Quantitative Finance, Higher School
of Economics,
Moscow, Russia.\\
\printead{e1}}
\end{aug}

\received{\smonth{1} \syear{2014}}
\revised{\smonth{6} \syear{2014}}

%
\begin{abstract}
We consider the problem of construction of goodness-of-fit tests for
diffusion processes with a \textit{small noise}. The basic hypothesis
is composite
parametric and our goal is to obtain asymptotically distribution-free
tests. We propose two solutions. The first one is based on a change of time,
and the second test is obtained using a linear transformation of
the ``natural'' statistics.
\end{abstract}

%
\end{frontmatter}

\section{Introduction}

We consider the following problem. Suppose that we observe a trajectory
$X^\varepsilon=\{X_t, 0\leq t\leq T\}$ of the following diffusion process:
%
%
\begin{equation}
\label{0} \mathrm{d}X_t=S (t,X_t ) \,\mathrm{d}t+
\varepsilon\sigma (t,X_t ) \,\mathrm{d}W_t,\qquad
X_0=x_0, 0\leq t\leq T,
\end{equation}
where $W_t,0\leq t\leq T$ is a Wiener process, $\sigma (t,x
)$ is
known smooth function, the initial value $x_0$ is deterministic
and the trend coefficient $S (t,x )$ is a unknown function. Here
$\varepsilon\in (0,1 )$ is a given parameter. We have to
test the
composite (parametric) hypothesis
%
%
\begin{equation}
\label{00} \mathscr{H}_0 \dvtx \mathrm{d}X_t=S (
\vartheta,t,X_t ) \,\mathrm{d}t+\varepsilon\sigma (t,X_t )
\,\mathrm{d}W_t,\qquad X_0=x_0, 0\leq t\leq
T
\end{equation}
against alternative $\mathscr{H}_1 \dvtx  \mbox{not } \mathscr{H}_0 $.
Here
$S (\vartheta,t,x )$ is a known smooth function of
$\vartheta$ and
$x$. The parameter $\vartheta\in\Theta$ is unknown and the set
$\Theta
\subset{ \RR}^d $ is open and bounded. Let us fix some value $\alpha
\in (0,1 )$ and consider the class of tests of asymptotic
($\varepsilon\rightarrow0$) size
$\alpha$:
\[
\mathcal{K}_\alpha= \bigl\{\bar\psi_\varepsilon\dvtx
\Ex_{\vartheta} \bar\psi_\varepsilon=\alpha+\mathrm{o} (1 ) \bigr\} \qquad
\mbox{for all } \vartheta\in\Theta.
\]
The test $\bar\psi_\varepsilon=\bar\psi_\varepsilon
(X^\varepsilon
 ) $ is the probability to reject the hypothesis $\mathscr{H}_0 $ and
$\Ex_{\vartheta}$ stands for the mathematical expectation under
hypothesis $\mathscr{H}_0 $.

Our goal is to find goodness-of-fit (GoF) tests which are
\textit{asymptotically distribution free} (ADF), that is, we look for
a test statistics
whose limit distributions under null hypothesis do not depend on the
underlying
model given by the functions $S (\vartheta,t,x )$, $\sigma
 (t,x )$ and the parameter~$\vartheta$. This work is a
continuation of
the study Kutoyants \cite{r9}, where an ADF test was
proposed in the case of
simple basic hypothesis.

The behaviour of stochastic systems governed by such equations (called
\textit{perturbed dynamical systems}) is well studied, see, for
example, Freidlin and Wentzell \cite{r3} and
the references therein. Estimation theory (parametric and non-parametric)
for such models of observations
is also well developped, see, for example,
Kutoyants \cite{r8} and Yoshida \cite{r17,r18}.

Let us remind the well-known basic results in this problem for the
i.i.d. model. We start with
the simple hypothesis. Suppose that we observe $n$ i.i.d.
r.v.'s $ (X_1,\ldots,X_n )=X^n$ with a continuous
distribution function
$F (x )$, and the basic hypothesis is
\[
\mathscr{H}_0\dvtx F (x )\equiv F_0 (x ),\qquad x\in
\RR.
\]
Then the Cram\'er--von Mises statistic is
\[
D _n=n\int \bigl[\hat F_n (x )-F_0 (x )
\bigr]^2\, \mathrm{d}F_0 (x ),\qquad \hat F_n
(x )=\frac{1}{n}\sum_{j=1}^{n}\one
_{ \{
X_j<x \}},
\]
where $\hat F_n (x )$ is the empirical distribution
function. Denote by $\mathcal{K}_\alpha$ the class of tests of asymptotic
($n\rightarrow\infty$) size $\alpha\in (0,1 )$, that is,
\[
\mathcal{K}_\alpha= \bigl\{\bar\psi\dvtx \Ex_0 \bar\psi=
\alpha+\mathrm{o} (1 ) \bigr\}.
\]

We have the convergence (under hypothesis $\mathscr{H}_0 $)
\[
B_n (x )=\sqrt{n} \bigl(\hat F_n (x )-F_0
(x ) \bigr)\Longrightarrow B \bigl( F_0 (x ) \bigr),
\]
where $B (\cdot )$ is a Brownian bridge process. Hence, it
can be
shown that
\[
D_n\Longrightarrow\delta\equiv\int_{0}^{1}B
(s )^2\,\mathrm{d}s,
\]
and the \textit{Cram\'er--von Mises} test
\[
\psi_n \bigl(X^n \bigr)=\one _{ \{D _n>c_\alpha \}}\in
\mathcal{K}_\alpha, \qquad\Pb \{\delta> c_\alpha \}=\alpha
\]
is
\textit{asymptotically distribution-free} (ADF).

The situation changes in the case of parametric basic hypothesis:
\[
\mathscr{H}_0 \dvtx F (x )=F (\vartheta ,x ),\qquad \vartheta\in
\Theta,
\]
where $\Theta= (\alpha,\beta )$. If we
introduce the similar statistic
\[
\hat D_n=n\int_{-\infty}^{\infty} \bigl[\hat
F_n (x )-F (\hat\vartheta_n,x ) \bigr]^2 \,
\mathrm{d}F (\hat\vartheta_n,x ),
\]
where $\hat\vartheta_n$ is the maximum likelihood estimator (MLE),
then (under
regularity conditions) we have
\[
U_n (x )=\sqrt{n} \bigl(\hat F_n (x )-F(\hat \vartheta
_n,x) \bigr)=B_n (x )-\sqrt{n}(\hat\vartheta_n-
\vartheta )\dot{F} (\vartheta,x )+\mathrm{o} (1 ).
\]
For the MLE, we can use its representation
\[
\sqrt{n} (\hat\vartheta_n-\vartheta ) =\frac{1}{\sqrt{n}}\sum
_{j=1}^ {n}\frac{\dot{\boldsymbol{\ell
}}
(\vartheta
,X_j )}{\mathrm{I} (\vartheta )} +\mathrm{o} (1
),\qquad \ell (\vartheta ,x )=\ln f (\vartheta ,x ).
\]
All this allows us to write the limit $U (\cdot )$ of the statistic
$U_n (\cdot )$
as follows:
\begin{eqnarray*}
U_n (x )&\Longrightarrow& B \bigl(F (\vartheta ,x ) \bigr)-\int
\frac{\dot
\ell (\vartheta,y )}{\sqrt{\mathrm{I} (\vartheta
 )}}\,\mathrm{d}B \bigl(F (\vartheta,y ) \bigr) \int
_{-\infty}^{x
}\frac{\dot
\ell (\vartheta
,y )}{ \sqrt{\mathrm{I} (\vartheta )}}\,\mathrm{d} F (
\vartheta,y )
\\
&=&B (t )-\int_{0}^{1}h (\vartheta,v )\,
\mathrm{d}B (v )\int_{0}^{t}h (\vartheta,v ) \,
\mathrm{d}v\equiv U (\vartheta,t ),
\end{eqnarray*}
where $t=F (\vartheta,x ) $ and we put
$ h (\vartheta,t )=\mathrm{I} (\vartheta
 )^{-1/2}{\dot\ell (\vartheta
, F^{-1}_\vartheta (t ) )}$.

If $\vartheta\in\Theta\subset\RR^d$, then we obtain a
similar equation
%
%
\begin{equation}
\label{iid} U (\vartheta,t )=B (t )- \biggl\langle\int_{0}^{1}
\mathbf{ h} (\vartheta,v ) \,\mathrm{d}B (v ),\int_{0}^{t}
\mathbf{ h} (\vartheta,v ) \,\mathrm{d}v \biggr\rangle,
\end{equation}
where $\langle\cdot,\cdot\rangle$ is the scalar product in $\RR^d$.

This presentation of the limit process $U (\vartheta,t )$
can be
found in Darling \cite{r2}. Of course, the test $\hat
\psi
_n=\one _{ \{\hat D_n>c_\alpha \}}$ is not ADF and the choice of
the threshold $c_\alpha$ can be a difficult problem. One\vspace*{2pt} way to avoid this
problem is, for example, to find a transformation $
L_W [U ] (t )=w (t )$, where $w
(\cdot
 )$ is the Wiener process. This transformation allows to write
the equality
\[
\Delta=\int_{-\infty}^{\infty}L_W [U ] \bigl(F
(\vartheta ,x ) \bigr)^2\,\mathrm{d} F (\vartheta ,x )=\int
_{0}^{1}w (t )^2\,\mathrm{d}t.
\]
Hence, if we prove the convergence
\[
\tilde D _n=\int_{-\infty}^{\infty
}L_W
[U_n ] (x )^2\,\mathrm{d}F (\hat \vartheta
_n,x )\Longrightarrow\Delta,
\]
then the test
$
\tilde\psi_n=\one _{ \{\tilde D_n>c_\alpha \}}$, with $
\Pb (\Delta>c_\alpha )=\alpha$
is ADF. Such transformation was proposed in Khmaladze \cite{r6}.

In the present work, we consider a similar problem for the model of
observations \eqref{0} with parametric basic hypothesis \eqref{00}.
Note that several problems of GoF testing for the model of observations
\eqref{0} with
simple basic hypothesis $\Theta= \{\vartheta_0 \}$ were
studied in
Dachian and Kutoyants \cite{r1}, Iacus and Kutoyants \cite{r5},
Kutoyants \cite{r9}. The tests considered there are mainly
based on the normalized difference $\varepsilon^{-1}
(X_t-x_t )$,
where $x_t=x_t (\vartheta_0 )$ is a solution of equation\vspace*{2pt}
\eqref{00}
for $\varepsilon=0$. This statistic is in some sense similar to the normalized
difference $\sqrt{n} (\hat F_n (x )-F_0 (x
) )$
used in the GoF problems for i.i.d. models. We propose two GoF
ADF tests. Note that the construction of the first test is in some
sense close
to the one considered in Kutoyants \cite{r12} and based on
the score
function process. These tests are originated by the different processes but
after our first transformation of the normalized difference
$\varepsilon
^{-1}(X_t-x_t(\hat\vartheta_\varepsilon))$ we obtain the same
integrals to
calculate as those in Kutoyants \cite{r12}.

Let us remind the related results in the case of simple hypothesis (see
Kutoyants \cite{r9}).
Suppose that the observed homogeneous diffusion process under null
hypothesis is
\[
\mathrm{d}X_t=S_0 (X_t ) \,\mathrm{d}t+
\varepsilon\sigma (X_t ) \,\mathrm{d}W_t,\qquad
X_0=x_0, 0\leq t\leq T,
\]
where $S_0 (x )$ is a known smooth function. Denote
$x_t= X_t\rrvert _{\varepsilon=0}$. We have $X_t\rightarrow x_t$ as
$\varepsilon\rightarrow0$ and we construct a GoF test based
on statistic $v_\varepsilon
 (t )=\varepsilon^{-1} (X_t-x_t )$. The limit of this
statistic is a Gaussian process. This process can be transformed into
the Wiener
process as
follows: introduce the statistic
\[
\delta_\varepsilon= \biggl[\int_{0}^{T}
\biggl(\frac{\sigma
 (x_t )}{S_0 (x_t )} \biggr)^{2}\,\mathrm{d}t \biggr]^{-2}
\int_{0}^{T} \biggl(\frac{X_t-x_t}{\varepsilon
S_0 (x_t )^{2}}
\biggr)^2 \sigma (x_t )^2 \,\mathrm{d}t.
\]
The following convergence:
\[
\delta_\varepsilon\Longrightarrow\Delta=\int_{0}^{1}w
(s )^2\,\mathrm{d}s
\]
was proved and therefore the test $\hat\psi_\varepsilon=\one _{ \{
\delta_\varepsilon>c_\alpha
 \}}$ with $
\Pb (\Delta>c_\alpha )=\alpha$ is ADF.

Consider now the hypotheses testing problem \eqref{0} and \eqref{00}. The\vspace*{2pt} solution
$x_t$ of equation \eqref{00} for $\varepsilon=0$ depends on
$\vartheta\in
\Theta\subset\RR^d$,
that is, $x_t=x_t (\vartheta )$. The statistic
$\hat v_\varepsilon (t )=\varepsilon^{-1}(X_t-x_t(\hat
\vartheta
_\varepsilon))$ (here $\hat\vartheta_\varepsilon$ is the MLE) is in some
sense similar to
$U_n (\cdot )$. Denote by $\hat v (t )$ the
limit of $\hat v_\varepsilon (t )$
as $\varepsilon\rightarrow0$ and suppose that we know the
transformation $L_U [\hat v  ]  (\cdot ) $ of
$\hat v (\cdot )$ into the Gaussian process
\[
U (\vartheta,t )=W (t )- \biggl\langle\int_{0}^{1}
\mathbf{h} (\vartheta,s ) \,\mathrm{d}W (s ),\int_{0}^{t}
\mathbf{h} (\vartheta,s ) \,\mathrm{d}s \biggr\rangle,\qquad 0\leq t\leq1
\]
with a vector-function $h (\vartheta,s )$ satisfying
\[
\int_{0}^{1}\mathbf{h} (\vartheta,s ) \mathbf{h}
(\vartheta ,s )^* \,\mathrm{d}s =\JJ.
\]
Here $\JJ$ is the $d\times d$ unit matrix.

The next steps are two transformations of $U (\cdot )$: one
transformation into
the Brownian bridge $L_B [U ] (s )=B
(s )$ and
another one into the Wiener process
$L_W [U ] (s )=w (s )$, respectively.
This allows
us to construct the ADF GoF tets as follows: let us introduce
(formally) the
statistics
\[
\delta_\varepsilon =\int_{0}^{T}
\bigl(L_B \bigl[L_U [\hat v_\varepsilon ] \bigr] (t )
\bigr)^2\,\mathrm{d}t,\qquad \Delta _\varepsilon =\int
_{0}^{T} \bigl(L_W
\bigl[L_U [\hat v_\varepsilon ] (t ) \bigr] \bigr)^2\,
\mathrm{d}t,
\]
and suppose that we have proved the convergences
\[
\delta_\varepsilon\Longrightarrow\delta =\int_{0}^{1}B
(s )^2\,\mathrm{d}s, \qquad\Delta_\varepsilon \Longrightarrow
\Delta=\int_{0}^{1}w (s )^2\,
\mathrm{d}s.
\]
Then the tests
\[
\hat\psi _\varepsilon=\one _{ \{ \delta
_\varepsilon>d_\alpha \}}, \qquad\Pb (\delta>d_\alpha
)=\alpha, \qquad\hat\Psi_\varepsilon=\one _{ \{ \Delta_\varepsilon
>c_\alpha
 \}}, \qquad\Pb (
\Delta>c_\alpha )=\alpha,
\]
belong to the class $\mathcal{K}_\alpha$ and are ADF. Our objective
is to realize
this program.

A similar result for ergodic diffusion processes is contained in
Kutoyants \cite{r10} (simple basic hypothesis) and
Kleptsyna and Kutoyants \cite{r7} (parametric basic
hypothesis).

\section{Auxiliary results}

We have the following stochastic differential equation:
%
%
\begin{equation}
\label{2} \mathrm{d}X_t=S (\vartheta,t,X_t ) \,
\mathrm{d}t+\varepsilon \sigma (t,X_t ) \,\mathrm{d}W_t,
\qquad X_0=x_0, 0\leq t\leq T,
\end{equation}
where $\vartheta\in\Theta$, $\Theta$ is an open bounded
subset of ${\RR}^d $ and $\varepsilon$ is a \textit{small parameter},
that is, we study this equation in the asymptotics of
\textit{small noise}
$\varepsilon\rightarrow0$.

Introduce the Lipschitz condition and that of linear growth:

$\mathcal{C}1$. \textit{The functions
$S (\vartheta,t,x )$ and $\sigma (t,x )$
satisfy the relations}
\begin{eqnarray*}
\bigl\llvert S (\vartheta,t,x )-S (\vartheta ,t,y )\bigr\rrvert +\bigl\llvert
\sigma (t,x )-\sigma (t,y )\bigr\rrvert &\leq& L\llvert x-y\rrvert ,
\\
\bigl\llvert S (\vartheta,t,x )\bigr\rrvert +\bigl\llvert \sigma (t,x )\bigr
\rrvert &\leq& L \bigl(1+\llvert x\rrvert \bigr).
\end{eqnarray*}
Recall that by these conditions the stochastic differential equation
\eqref{2}
has a unique strong solution (Liptser and Shiryaev \cite{r14}), and
moreover this solution
$X^\varepsilon= \{X_t,0\leq t\leq T \}$ converges
uniformly, with
respect to $t$, to the solution $x^T= \{x_t,0\leq t\leq T \}
$ of the
ordinary differential equation
%
%
\begin{equation}
\label{3} \frac{\mathrm{d}x_t}{\mathrm{d}t}=S (\vartheta,t,x_t ),\qquad
x_0, 0\leq t\leq T.
\end{equation}
Observe that $x_t=x_t (\vartheta )$ (for the proof see
Freidlin and Wentzell \cite{r3}, Kutoyants \cite{r8}).

$\mathcal{C}2$. \textit{The diffusion coefficient $\sigma
(t,x )^2$
is bounded away from zero}
\[
\inf_{0\leq t\leq T,x}\sigma (t,x )^2>0.
\]
Conditions $\mathcal{C}1$ and $\mathcal{C}2 $ provide the equivalence of
the measures
$ \{\Pb_\vartheta^{ (\varepsilon )},\vartheta\in
\Theta \}$ induced on the
measurable space $ (\mathscr{C}_T,{\goth
B}_T  )$ by the solutions of equation \eqref{2}
(Liptser and Shiryaev \cite{r14}). Here $\mathscr{C}_T $ is the space of
continuous functions on
$ [0,T ]$
with uniform metrics and ${\goth
B}_T $ is the Borelian $\sigma$-algebra of its subsets. The
likelihood ratio function is
\[
L \bigl(\vartheta,X^\varepsilon \bigr)=\exp \biggl\{\int_{0}^{T}
\frac{S (\vartheta
,t,X_t ) }{\varepsilon^2\sigma
 (t,X_t )^2} \,\mathrm{d}X_t-\int_{0}^{T}
\frac{S
(\vartheta
,t,X_t )^2 }{2\varepsilon^2\sigma
 (t,X_t )^2} \,\mathrm{d}t \biggr\},\qquad \vartheta\in\Theta,
\]
and the maximum likelihood estimator (MLE) $\hat\vartheta_\varepsilon
$ is
defined by the equation
\[
L \bigl(\hat\vartheta_\varepsilon,X^\varepsilon \bigr)=\sup
_{\vartheta
\in\Theta}L \bigl(\vartheta,X^\varepsilon \bigr).
\]

The following regularity conditions (smoothness and identifiability) provides
us necessary properties of the MLE.
Below $x_t=x_t (\vartheta_0 )$.

$\mathcal{C}3$. \textit{The functions $S (\vartheta,t,x
)$ and $\sigma
 (t,x )$ have two continuous bounded derivatives w.r.t. $x$
and the
function $S (\vartheta,t,x )$ has two continuous bounded derivatives
w.r.t. $\vartheta$}.

\textit{For any $\nu>0$}
\[
\inf_{\vartheta_0 \in\Theta}\inf_{\llvert \vartheta
-\vartheta_0\rrvert  >\nu}\int
_{0}^{T} \biggl(\frac{S
(\vartheta
,t,x_t )-S (\vartheta_0
,t,x_t ) }{\sigma
 (t,x_t )}
\biggr)^2\,\mathrm{d}t >0
\]
\textit{and the information matrix ($d\times d$)
\[
{\II} (\vartheta_0 )=\int_{0}^{T}
\frac{\dot\mathbf{S} (\vartheta_0
,t,x_t ) \dot\mathbf{S} (\vartheta_0, t,x_t
)^*}{\sigma
 (t,x_t )^2} \,\mathrm{d}t
\]
is uniformly non-degenerate:}
\[
\inf_{\vartheta_0 \in\Theta}\inf_{\llvert \lambda\rrvert =1} \lambda^* \II (
\vartheta_0 )\lambda>0.
\]

We denote by a prime the derivatives w.r.t. $x$ and $t$, and by a
dot those w.r.t. $\vartheta$, that is, for a function
$f=f (\vartheta,t,x )$ we write
\begin{eqnarray*}
f' (\vartheta,t,x )&=&\frac{\partial f (\vartheta
,t,x )}{\partial x},
\\
f_t' (\vartheta,t,x )&=&\frac
{\partial
f (\vartheta,t,x )}{\partial t},
\\
\dot f (\vartheta ,t,x )&=&\frac{\partial
f (\vartheta,t,x )}{\partial\vartheta}.
\end{eqnarray*}
Of course, in the case of $d>1 $ the
derivative $\dot\mathbf{f} (\vartheta,t,x )$ is a column vector.

If the conditions $\mathcal{C}2 $ and $\mathcal{C}3$
hold, then the MLE admits the representation
%
%
\begin{eqnarray}
\label{mle} \varepsilon^{-1} ({\hat\vartheta_\varepsilon-
\vartheta}  )=\II (\vartheta )^{-1} \int_{0}^{T}
\frac{\dot\mathbf{S} (\vartheta
,t,x_t ) }{\sigma
 (t,x_t )} \,\mathrm{d}W_t +\mathrm{o} (1 ).
\end{eqnarray}
Here, $x_t=x_t (\vartheta )$. For the proof see,
Kutoyants \cite{r9}.

Note that $X_t=X_t (\varepsilon )$ (solution of equation
\eqref{2}) under condition $\mathcal{C}3$ is continuously differentiable
w.r.t. $\varepsilon$. Denote the derivatives
\[
X^{ (1 )}_t=\frac{\partial X_t}{\partial\varepsilon},\qquad x^{ (1 )}_t=
 \frac{\partial X_t}{\partial\varepsilon
}\Big\rrvert _{\varepsilon=0},\qquad 0\leq t\leq T.
\]
The equations for $X_t^{ (1 )}$ and $x_t^{ (1 )}$ are
\[
\mathrm{d}X_t^{ (1 )}=S' (\vartheta
,t,X_t )X_t^{ (1 )} \,\mathrm{d}t+ \bigl[
\varepsilon \sigma' (t,X_t )X_t^{ (1 )}+
\sigma (t,X_t ) \bigr] \,\mathrm{d}W_t,\qquad
X_0^{ (1 )}=0
\]
and
%
%
\begin{equation}
\label{4} \mathrm{d}x_t^{ (1 )}=S' (
\vartheta ,t,x_t )x_t^{ (1 )} \,\mathrm{d}t+\sigma
(t,x_t ) \,\mathrm{d}W_t,\qquad x_0^{ (1 )}=0,
\end{equation}
respectively.
Hence $x_t^{ (1 )}$, $0\leq t\leq T $ is a Gaussian process and
it can
be written as
%
%
\begin{equation}
\label{x1} x_t^{ (1 )}=\int_{0}^{t}
\exp \biggl\{\int_{s}^{t}S' (
\vartheta ,v,x_v )\,\mathrm{d}v \biggr\} \sigma (s,x_s )
\,\mathrm{d}W_s.
\end{equation}

Denote
\[
\psi (t )=\exp \biggl\{\int_{0}^{t}S'
(\vartheta ,v,x_v )\,\mathrm{d}v \biggr\},\qquad \psi_\varepsilon
(t )=\exp \biggl\{\int_{0}^{t}S' (
\hat\vartheta _\varepsilon,v,X_v )\,\mathrm{d}v \biggr\}.
\]

We can write
\begin{eqnarray*}
\frac{X_t-x_t(\hat\vartheta_\varepsilon)}{\varepsilon
}&=&\frac{X_t-x_t(\vartheta)}{\varepsilon}+\frac{x_t (\vartheta
 )-x_t (\hat\vartheta_\varepsilon )}{\varepsilon}
\\
&=&X_t^{ (1 )}- \biggl\langle\frac{(\hat\vartheta
_\varepsilon-\vartheta
)}{\varepsilon}, \dot
\mathbf{x}_t (\vartheta ) \biggr\rangle+\mathrm{o} (1 )
\\
&=&x_t^{ (1 )}- \biggl\langle \II (\vartheta )^{-1}
\int_{0}^{T}\frac{\dot\mathbf{S} (\vartheta
,s,x_s ) }{\sigma (s,x_s )} \,
\mathrm{d}W_s, \dot \mathbf{x}_t (\vartheta ) \biggr
\rangle+\mathrm{o} (1 )
\\
&=&\psi (t )V (t )+\mathrm{o} (1 ),
\end{eqnarray*}
where
\[
V (t )= \psi (t )^{-1} x_t^{ (1 )}-\psi (t
)^{-1} \biggl\langle\II (\vartheta )^{-1} \int
_{0}^{T}\frac{\dot\mathbf{S} (\vartheta
,s,x_s ) }{\sigma
 (s,x_s )} \,\mathrm{d}W_s,
\dot\mathbf{x}_t (\vartheta ) \biggr\rangle.
\]
Introduce the random process
\[
U (\vartheta,t )= \int_{0}^{t}\frac{\psi (s )
}{\sigma
 (s,x_s )}
\,\mathrm{d}V (s ).
\]

%
\begin{lemma}\label{L1}
We have the equality
%
%
\begin{equation}
\label{8} U (\vartheta,t )=W_t- \biggl\langle\int
_{0}^{T}\mathbf{h} (\vartheta ,s )\,
\mathrm{d}W_s,\int_{0}^{t}\mathbf{h} (
\vartheta ,s ) \,\mathrm{d}s \biggr\rangle,\qquad 0\leq t\leq T,
\end{equation}
where
%
%
\begin{equation}
\label{h} \mathbf{h} (\vartheta,t )=\II (\vartheta )^{-1/2}
\frac{\dot
\mathbf{S} (\vartheta,t,x_t )}{\sigma (t,x_t )}
\end{equation}
is a vector-valued function.
\end{lemma}

\begin{pf}
The solution of equation \eqref{4} can be written (see \eqref{x1}) as
\[
x_t^{ (1 )}=\int_{0}^{t}
\frac{\psi (t )\sigma
 (s,x_s )}{\psi
 (s )}\,\mathrm{d}W_s.
\]
For the vector $\dot\mathbf{x}_t (\vartheta )$, we can write
\[
\dot\mathbf{x}_t (\vartheta )=\int_{0}^{t}S'
(\vartheta ,s,x_s ) \dot\mathbf{x}_s (\vartheta ) \,
\mathrm{d}s + \int_{0}^{t}\dot\mathbf{S} (
\vartheta ,s,x_s )\,\mathrm{d}s.
\]
The solution of this equation is
\[
\dot\mathbf{x}_t (\vartheta )=\psi (t )\int_{0}^{t}
\frac{\dot
\mathbf{S} (\vartheta
,s,x_s )}{\psi (s )} \,\mathrm{d}s.
\]
Introduce two stochastic processes
\[
v_1 (t )=\psi (t )^{-1}x_t^{ (1
)}=
\int_{0}^{t} \psi (s )^{-1} \sigma
(s,x_s ) \,\mathrm{d}W_s
\]
and
\[
\mathbf{v}_2 (t )=\psi (t )^{-1}\dot\mathbf{x}_t
(\vartheta )=\int_{0}^{t}\psi (s )^{-1}
\dot\mathbf{S} (\vartheta ,s,x_s ) \,\mathrm{d}s.
\]
Then we can write
\begin{eqnarray*}
U (\vartheta,t )&=&\int_{0}^{t}
\frac{\psi (s )
}{\sigma (s,x_s )} \,\mathrm{d} V (s )
\\
&=&\int_{0}^{t}\frac{\psi (s ) }{\sigma
 (s,x_s )} \,\mathrm{d}
v_1 (s )
\\
&&{} - \biggl\langle \II (\vartheta )^{-1} \int_{0}^{T}
\frac{\dot\mathbf{S} (\vartheta
,s,x_s ) }{\sigma (s,x_s )} \,\mathrm{d}W_s, \int_{0}^{t}
\frac{\psi (s ) }{\sigma (s,x_s
)} \,\mathrm{d} \mathbf{v}_2 (s ) \biggr\rangle
\\
&=&W (t )- \biggl\langle \II (\vartheta )^{-1/2} \int
_{0}^{T}\frac{\dot\mathbf{S} (\vartheta
,s,x_s ) }{\sigma (s,x_s )} \,\mathrm{d}W_s,
\II (\vartheta )^{-1/2}\int_{0}^{t}
\frac{\dot\mathbf{S} (\vartheta
,s,x_s ) }{\sigma (s,x_s )} \,\mathrm{d} s \biggr\rangle
\\
&=&W_t- \biggl\langle\int_{0}^{T}
\mathbf{h} (\vartheta,s )\,\mathrm{d}W_s,\int_{0}^{t}
\mathbf{h} (\vartheta,s ) \,\mathrm{d}s \biggr\rangle.
\end{eqnarray*}

Introduce the random process
\[
u (\vartheta,r )=T^{-1/2}U (\vartheta,rT ),\qquad 0\leq r\leq1
\]
and denote
\begin{eqnarray*}
\II_1 (\vartheta )&=&\int_{0}^{1}
\frac{\dot\mathbf{S} (\vartheta
,rT,x_{rT} )\dot\mathbf{S} (\vartheta,rT,x_{rT}
)^*}{\sigma
 (rT,x_{rT} )^2} \,\mathrm{d}r,
\\
\tilde\mathbf{h} (\vartheta,r )&=&\II_1 (\vartheta )^{-1/2}
\frac{\dot\mathbf{S} (\vartheta
,rT,x_{rT} )}{\sigma (rT,x_{rT} )},\qquad w_r=T^{-1/2}W_{rT}.
\end{eqnarray*}

Then we can write
%
%
\begin{equation}
\label{tu} u (\vartheta,r )=w_r- \biggl\langle\int
_{0}^{1}\tilde \mathbf{h} (\vartheta,q )
\,\mathrm{d}w_q,\int_{0}^{r}\tilde
\mathbf{h} (\vartheta,q ) \, \mathrm{d}q \biggr\rangle,\qquad 0\leq r\leq1,
\end{equation}
and therefore
\[
\int_{0}^{1 }\tilde\mathbf{h} (\vartheta,q )
\tilde \mathbf{h} (\vartheta,q )^* \,\mathrm{d}q =\JJ.
\]
\upqed\end{pf}

Note that $u (\cdot )$ is in some sense a \textit{universal limit} which
appears in the problems of goodness of fit testing for stochastic
processes. For example, the same limit is obtained in the case of ergodic
diffusion process and in the case of inhomogeneous Poisson process
(Kutoyants \cite{r11}). The main difference with the i.i.d.
case is due to the Wiener
process here, while in the i.i.d. case the Brownian bridge $B(t)$, $0\leq
t\leq1$ appears (see \eqref{iid}). Of course, we can immediately replace
$B (t )$ by a Wiener
process $B (t )=W_t-W_1 t$ and this will increase the
dimension of
the vector $h (\vartheta,\cdot )$. In the case of
vector-valued parameter
$\vartheta$, this change is not essential and will slightly modify
calculations of the test statistics for the first type test. At the
same time if
the parameter $\vartheta$ is one-dimensional, then we can easily construct
the second-type goodness-of-fit test for stochastic processes and it remains
unclear how to construct such tests in the i.i.d. case. The difference
will be
explained in Section~\ref{sec3.2}.

In the construction of a GoF test, we will use another condition.

$\mathcal{C}4$. \textit{The functions $S (\vartheta,t,x
 )$, $\dot\mathbf{S} (\vartheta,t,x
 )$ and $\sigma (t,x
 )$ have continuous bounded derivatives w.r.t. $t\in
[0,T ]$}.

\section{Main results}

Suppose that we observe a trajectory $X^\varepsilon
= (X_t,0\leq t\leq T )$ of the following diffusion process:
%
%
\begin{equation}
\label{5} \mathrm{d}X_t=S (t,X_t ) \,\mathrm{d}t+
\varepsilon\sigma (t,X_t ) \,\mathrm{d}W_t,\qquad
X_0=x_0, 0\leq t\leq T.
\end{equation}

We have to test the basic parametric hypothesis
\[
\mathcal{H}_0\dvtx S (t,x )=S (\vartheta ,t,x ),\qquad 0\leq t\leq
T, \vartheta\in\Theta,
\]
that is, that the observed process
\eqref{5} has the stochastic differential
%
%
\begin{equation}
\label{6} \mathrm{d}X_t=S (\vartheta,t,X_t ) \,
\mathrm{d}t+\varepsilon \sigma (t,X_t ) \,\mathrm{d}W_t,
\qquad X_0=x_0, 0\leq t\leq T
\end{equation}
with some $\vartheta\in\Theta$. Here $S (\vartheta,t,x
)$ and
$\sigma (t,x )$ are known strictly positive smooth
functions and
$\Theta\subset R^d$ is an open convex set. We have to test this
hypothesis in
the asymptotics of a \textit{small noise} (as $\varepsilon\rightarrow0$).

Our goal is to construct such statistics $v_\varepsilon[X^\varepsilon
] (\cdot ) $, $V_\varepsilon[X^\varepsilon
] (\cdot ) $ that (under hypothesis $\mathcal{H}_0 $)
\begin{eqnarray*}
\delta_\varepsilon&=&\int_{0}^{T}v_\varepsilon
\bigl[X^\varepsilon \bigr] (t )^2\,\mathrm{d}t\quad\Longrightarrow
\quad\delta=\int_{0}^{1}B (s )^2\,
\mathrm{d}s,
\\
\Delta_\varepsilon&=&\int_{0}^{T}V_\varepsilon
\bigl[X^\varepsilon \bigr] (t )^2\,\mathrm{d}t\quad\Longrightarrow
\quad\Delta=\int_{0}^{1}w (s )^2\,
\mathrm{d}s,
\end{eqnarray*}
where $B (\cdot )$ and $w (\cdot ) $ are the
Brownian bridge
and the Wiener process, respectively.
Then we introduce the tests
\[
\hat\psi_\varepsilon=\one _{ \{\delta_\varepsilon>d_\alpha
 \}},\qquad \hat\Psi_\varepsilon=\one
_{ \{\Delta_\varepsilon
>c_\alpha
 \}}
\]
with the thresholds $c_\alpha$ and $d_\alpha$ satisfying the equations
%
%
\begin{equation}
\label{7} \Pb (\delta>d_\alpha )=\alpha,\qquad \Pb (\Delta
>c_\alpha )=\alpha.
\end{equation}
These tests will belong to the class
\[
\mathcal{K}_\alpha= \Bigl\{\bar\psi_\varepsilon\dvtx \lim
_{\varepsilon
\rightarrow0}\Ex_\vartheta\bar\psi_\varepsilon=\alpha, \forall
\vartheta\in\Theta \Bigr\}
\]
and will be ADF.

We propose these tests in the Sections~\ref{sec3.1} and~\ref{sec3.2} below. We call $\hat
\psi_\varepsilon$
\textit{the first} test and $\hat\Psi_\varepsilon$ \textit{the
second} test.

\subsection{First test}\label{sec3.1}

The construction of the first ADF GoF test is based on the following well
known
property. Suppose that we have a Gaussian process $U (t ),
0\leq
t\leq T$ satisfying the equation
\[
U (t )=w (t )-\int_{0}^{T}h (s ) \,\mathrm{d}w
(s ) \int_{0}^{t}h (s )\,\mathrm{d}s,\qquad \int
_{0}^{T}h (s )^2 \,\mathrm{d}s=1.
\]

Introduce the process
\begin{eqnarray*}
b (t )&=&\int_{0}^{t}h (s ) \,\mathrm{d}U (s )
\\
&=&\int_{0}^{t}h (s ) \,\mathrm{d}w (s )-\int
_{0}^{T}h (s ) \,\mathrm{d}w (s )\int
_{0}^{t}h (s )^2 \,\mathrm{d}s.
\end{eqnarray*}
It is easy to see that $b (0 )=b (T )=0$ and
\[
\Ex \bigl[b (t )b (s ) \bigr]=\int_{0}^{t\wedge
s}h (v
)^2 \,\mathrm{d}v-\int_{0}^{t}h (v
)^2 \,\mathrm{d}v\int_{0}^{ s}h (v
)^2 \,\mathrm{d}v.
\]

Let us put
\[
\tau=\int_{0}^{ s}h (v )^2 \,
\mathrm{d}v,\qquad b (t )=B (\tau ),\qquad 0\leq\tau\leq1.
\]
Then
\begin{eqnarray*}
\delta&=&\int_{0}^{T} \biggl(\int
_{0}^{t}h (s ) \,\mathrm{d}U (s )
\biggr)^2h (t )^2 \,\mathrm{d}t
\\
&=&\int_{0}^{T}b (t )^2h (t
)^2 \,\mathrm{d}t =\int_{0}^{1}B (
\tau )^2 \,\mathrm{d}\tau.
\end{eqnarray*}

Suppose that the parameter $\vartheta$ is one-dimensional, $\vartheta
\in\Theta
= (a,b )$ and that we already proved the convergence (see
Lemma~\ref{L1})
\[
U_\varepsilon (t )=\int_{0}^{t}
\frac{\psi_\varepsilon
 (s )}{\sigma (s,X_s )} \,\mathrm{d} \biggl(\frac{X_s-x_s(\hat\vartheta_\varepsilon)}{\varepsilon\psi
_\varepsilon (s )} \biggr)\longrightarrow U (
\vartheta ,t ),\qquad 0\leq t\leq T,
\]
where
\[
U (\vartheta,t )=w (t )-\int_{0}^{T}h (
\vartheta,s ) \,\mathrm{d}w (s ) \int_{0}^{t}h (
\vartheta,s )\,\mathrm{d}s,\qquad \int_{0}^{T}h (
\vartheta,s )^2 \,\mathrm{d}s=1.
\]
Recall that
\[
h (\vartheta,s )=\mathrm{I} (\vartheta )^{-1/2}\frac{\dot
S (\vartheta,s,x_s )}{ \sigma (s,x_s )},\qquad
\mathrm{I} (\vartheta )=\int_{0}^{T}
\frac{\dot
S (\vartheta,s,x_s )^2}{ \sigma (s,x_s
)^2}\,\mathrm{d}s.
\]

Introduce (formally) the statistic
\[
\hat\delta_{\varepsilon}=\int_{0}^{T} \biggl(
\int_{0}^{t}h (\hat\vartheta_\varepsilon,s ) \,
\mathrm{d}U_{\varepsilon} (s ) \biggr)^2h(\hat\vartheta
_\varepsilon ,t)^2 \,\mathrm{d}t.
\]
If we prove that
\begin{eqnarray*}
&&\int_{0}^{T} \biggl(\int_{0}^{t}h
(\hat\vartheta_\varepsilon,s ) \,\mathrm{d}U_{\varepsilon} (s )
\biggr)^2h(\hat\vartheta _\varepsilon ,t)^2 \,
\mathrm{d}t
\\
&&\quad \Longrightarrow\int_{0}^{T} \biggl(\int
_{0}^{t}h (\vartheta,s ) \,\mathrm{d}U (
\vartheta,s ) \biggr)^2h(\vartheta ,t)^2 \,\mathrm{d}t
\end{eqnarray*}
then the test $\hat\psi_\varepsilon=\one _{ \{\delta_\varepsilon
>c_\alpha \}}$
will be {ADF}.

The main technical problem in carrying out this program is to
define the stochastic integral
\[
\int_{0}^{t}h (\hat\vartheta_\varepsilon,s )
\,\mathrm{d}U_{\varepsilon} (s )
\]
containing the MLE
$\hat\vartheta_\varepsilon=\hat\vartheta_\varepsilon
(X_t,0\leq t\leq
T )$. We will proceed as follows: First, we formally
differentiate and integrate and then we take the final expressions, which
do not contain stochastic integrals, as starting statistics.

Introduce the statistics
\begin{eqnarray*}
D(\vartheta,s,X_s)&=&S \bigl(\vartheta,s,x_s (\vartheta
) \bigr)+S' (\vartheta,s,X_s ) \bigl(X_s-x_s
(\vartheta ) \bigr),
\\
R \bigl(\vartheta,t,X^t \bigr)&=&\int_{x_0}^{X_t}
\frac{\dot
S(\vartheta,t,y)}{
\sqrt{\mathrm{I}(\vartheta)} \sigma
 (t,y )^2} \,\mathrm{d}y
\\
&&{} -\int_{0}^{t}\int_{x_0}^{X_s}
\frac{ \dot S'_s(\vartheta
,s,y)\sigma
 (s,y )-2\dot S(\vartheta,s,y)\sigma_s'
 (s,y )}{ \sqrt{\mathrm{I}(\vartheta)} \sigma
 (s,y )^3} \,\mathrm{d}y \,\mathrm{d}s,
\\
Q \bigl(\vartheta,t,X^t \bigr)&=&\int_{0}^{t}
\frac{\dot S(\vartheta
,s,X_s)D(\vartheta,s,X_s)}{\sqrt{\mathrm{I}(\vartheta)} \sigma
 (s,X_s )^2} \,\mathrm{d}s,
\\
K_\varepsilon ( \vartheta,t )& =&\varepsilon^{-1} \bigl[R \bigl(
\vartheta ,t,X^t \bigr)-Q \bigl(\vartheta,t,X^t \bigr)
\bigr],
\\
\delta_\varepsilon&=&\int_{0}^{T}K_\varepsilon(
\hat\vartheta _\varepsilon ,t)^2 h_\varepsilon(\hat
\vartheta_\varepsilon,t )^2 \,\mathrm{d}t.
\end{eqnarray*}
The first test is given in the following theorem.

%
\begin{theorem}
\label{T1}
Suppose that the conditions $\mathcal{C}1$--\hspace*{0.5pt}$\mathcal{C}4$ hold. Then
the test
\[
\hat\psi_\varepsilon=\one _{ \{\delta_\varepsilon>c_\alpha
 \}},\qquad \Pb \{\delta>c_\alpha
\} =\alpha
\]
is ADF and belongs to $\mathcal{K}_\varepsilon$.
\end{theorem}

\begin{pf}
We can write (formally)
%
%
\begin{eqnarray}
\label{u} U_\varepsilon (t )&=&\int_{0}^{t}
\frac{\psi_\varepsilon
 (s )}{\sigma (s,X_s )} \,\mathrm{d}V_\varepsilon (s )
\nonumber
\\
&=&\int_{0}^{t}\frac{\psi_\varepsilon (s
)}{\sigma
 (s,X_s )}\,\mathrm{d}
\biggl(\frac{X_s-x_s(\hat\vartheta
_\varepsilon
)}{\psi_\varepsilon
 (s )\varepsilon} \biggr)
\nonumber
\\[-8pt]
\\[-8pt]
& =&\int_{0}^{t}\frac{\mathrm{d}X_s}{\varepsilon\sigma
 (s,X_s )}- \int
_{0}^{t} \biggl[\frac{S(\hat\vartheta
_\varepsilon
,s,x_s(\hat\vartheta_\varepsilon))}{\varepsilon\sigma
 (s,X_s )}+
\frac{S'(\hat\vartheta_\varepsilon
,s,X_s) (X_s-x_s(\hat\vartheta_\varepsilon
) )}{\varepsilon\sigma
 (s,X_s )} \biggr] \,\mathrm{d}s
\nonumber
\\
& =&\int_{0}^{t}\frac{\mathrm{d}X_s}{\varepsilon\sigma
 (s,X_s )}- \int
_{0}^{t}\frac{D(\hat\vartheta
_\varepsilon
,s,X_s)}{\varepsilon\sigma
 (s,X_s ) }\,\mathrm{d}s,
\nonumber
\end{eqnarray}
where we have used the equality
\[
\mathrm{d}x_s(\hat\vartheta_\varepsilon)=S\bigl(\hat\vartheta
_\varepsilon,s, x_s(\hat\vartheta_\varepsilon)\bigr)\,
\mathrm{d}s.
\]

Hence (formally), we obtain the following expression.
\begin{eqnarray*}
\int_{0}^{t} h_\varepsilon(\hat
\vartheta_\varepsilon,s)\,\mathrm{d}U_\varepsilon (s ) &=&\int
_{0}^{t} \frac{\dot S(\hat
\vartheta
_\varepsilon,s,X_s)}{\sqrt{ \mathrm{I}(\hat\vartheta_\varepsilon)}
\varepsilon\sigma
 (s,X_s )^2} \,
\mathrm{d}X_s
\\
&&{} -\int_{0}^{t} \frac{\dot S(\hat\vartheta
_\varepsilon
,s,X_s)D(\hat\vartheta_\varepsilon,s,X_s)}{\sqrt{\mathrm{I}(\hat
\vartheta
_\varepsilon) } \varepsilon\sigma
 (s,X_s )^2} \,
\mathrm{d}s.
\end{eqnarray*}

The estimator $\hat\vartheta_\varepsilon=\hat\vartheta
_\varepsilon (X_t,0\leq t\leq T )$ and therefore the stochastic
integral is not well\vspace*{2pt} defined because the integrand $\dot S(\hat
\vartheta
_\varepsilon
,s,X_s) $ is not a non-anticipative random function. Note that in the
linear case $S (\vartheta,t,x )=\vartheta Q (s,x
)$ we have
no such problem (see example below). This difficulty can be avoided in general
case by at least two ways: The
first one is to replace the stochastic integral by
it's \textit{robust
version} as we show below. The second possibility is to use a consistent
estimator $\bar\vartheta_{\nu_\varepsilon}$ of the parameter
$\vartheta$
constructed after the observations
$X^{\nu_\varepsilon}= (X_t,0\leq t\leq\nu_\varepsilon
)$, where
$\nu_\varepsilon\rightarrow0$ but sufficiently slowly. With this
estimator, we can calculate the integral
\[
\int_{\nu_\varepsilon}^{t}\frac{\dot S (\bar\vartheta_{\nu
_\varepsilon
},s,X_s )}{\sigma (s,X_s )^2} \,
\mathrm{d}X_s
\]
without any problem, and all limits will be the same. Such construction is
discussed for a different problem in Kutoyants and Zhou \cite{r13}.

Introduce the function
\[
M (\vartheta,t,x )=\int_{x_0}^{x}
\frac{\dot S(\vartheta
,t,y)}{ \sigma
 (t,y )^2} \,\mathrm{d}y.
\]
Then by the It\^o formula
\begin{eqnarray*}
\mathrm{d}M(\vartheta,t,X_t)&=&M_t' (
\vartheta,t,X_t ) \,\mathrm{d}t+\frac{\varepsilon^2\sigma
 (t,X_t )^2}{2}
M_{xx}'' (\vartheta,t,X_t ) \,
\mathrm{d}t
\\
&&{} +M_x' (\vartheta,t,X_t )\,
\mathrm{d}X_t
\end{eqnarray*}
and therefore
\begin{eqnarray*}
&&\int_{0}^{t} \frac{\dot S(\vartheta,s,X_s)}{\sigma
 (s,X_s )^2} \,
\mathrm{d}X_s
\\
&&\quad=M(\vartheta ,t,X_t) -\int_{0}^{t}
\biggl[M_s' (\vartheta,s,X_s )+
\frac
{\varepsilon^2\sigma
 (s,X_s )^2}{2} M_{xx}'' (
\vartheta,s,X_s ) \biggr]\,\mathrm{d}s
\\
&&\quad =\int_{x_0}^{X_t}\frac{\dot S(\vartheta,t,y)}{ \sigma
 (t,y )^2} \,
\mathrm{d}y-\int_{0}^{t}\int_{x_0}^{X_s}
\frac
{\dot S_s'(\vartheta,s,y)}{ \sigma
 (s,y )^2} \,\mathrm{d}y
\\
&&\qquad{} + \int_{0}^{t}\int
_{x_0}^{X_s} \frac{2\dot S(\vartheta,s,y)\sigma_s'
 (s,y )}{ \sigma
 (s,y )^3} \,\mathrm{d}s-
\frac{\varepsilon^2}{2}\int_{0}^{t}{\sigma
(s,X_s )^2}   M_{xx}''
(\vartheta,s,X_s )\,\mathrm{d}s.
\end{eqnarray*}
Note that the contribution of the term
\[
\varepsilon^2 \int_{0}^{t} \sigma
(s,X_s )^2 M_{xx}'' (
\hat\vartheta_\varepsilon ,s,X_s ) \,\mathrm{d}s
\]
is asymptotically ($\varepsilon\rightarrow0$) negligible. Therefore,
\[
K_\varepsilon (\vartheta,t ) =\varepsilon^{-1} \bigl[R \bigl(
\vartheta ,t,X^t \bigr)-Q \bigl(\vartheta,t,X^t \bigr)
\bigr]
\]
is asymptotically equivalent to
\[
\tilde K_\varepsilon ( \vartheta,t )=\int_{0}^{t}
h_\varepsilon(\vartheta,s) \,\mathrm{d}U_\varepsilon (s ).
\]
The difference is in the dropped term of order $\mathrm{O} (\varepsilon
 )$.

We have to verify the convergence of the integrals
\[
\delta_\varepsilon=\int_{0}^{T}
\frac{K_\varepsilon(\hat\vartheta
_\varepsilon,t)^2 \dot
S(\hat\vartheta_\varepsilon,t,X_t)^2}{{\mathrm{I}(\hat\vartheta
_\varepsilon
)} \sigma (t,X_t )^2} \,\mathrm{d}t \longrightarrow \int_{0}^{T}
\frac{K(\vartheta,t)^2 \dot
S(\vartheta,t,x_t)^2}{{\mathrm{I}(\vartheta
)} \sigma (t,x_t )^2} \,\mathrm{d}t.
\]

Regularity conditions $\mathcal{C}1$--\hspace*{0.5pt}$\mathcal{C}3$ give the uniform
convergences
\begin{eqnarray*}
\sup_{0\leq t\leq T}\bigl\llvert X_t-x_t (
\vartheta )\bigr\rrvert &\longrightarrow&0,\qquad \mathrm{I}(\hat\vartheta
_\varepsilon) \longrightarrow\mathrm{I} (\vartheta ),
\\
\sup_{0\leq t\leq T}\bigl\llvert h_\varepsilon(\hat\vartheta
_\varepsilon,t)- h(\hat\vartheta_\varepsilon,t)\bigr\rrvert &=&\sup
_{0\leq t\leq T}\biggl\llvert \frac{\dot
S(\hat\vartheta_\varepsilon,t,X_t)}{\sqrt{\mathrm{I}(\hat\vartheta
_\varepsilon)}\sigma (t,X_t )}-\frac{\dot
S(\vartheta,t,x_t)}{\sqrt{\mathrm{I}(\vartheta)}\sigma
(t,x_t )}\biggr
\rrvert \longrightarrow0.
\end{eqnarray*}
Introduce two processes
\begin{eqnarray*}
Y_\varepsilon \bigl(\hat\vartheta_\varepsilon,t,X^t \bigr)&=&
\int_{0}^{t}\frac{\dot S (\hat\vartheta_\varepsilon
,s,X_s ) [ S (\vartheta,s,X_s ) -D(\hat
\vartheta
_\varepsilon,s,X_s)  ]}{\sigma (s,X_s )^2} \,\mathrm{d}s,
\\
Z \bigl(\hat\vartheta_\varepsilon,t,X^t \bigr)&=&{R \bigl(\hat
\vartheta _\varepsilon,t,X^t \bigr)}- \int_{0}^{t}
\frac{\dot S (\hat
\vartheta
_\varepsilon,s,X_s )S (\vartheta,s,X_s )}{\sigma
 (s,X_s )^2} \,\mathrm{d}s.
\end{eqnarray*}
Then
\[
K_\varepsilon (t )=\varepsilon^{-1} \bigl[Y_\varepsilon \bigl(
\hat\vartheta_\varepsilon,t,X^t \bigr)+ Z \bigl(\hat
\vartheta_\varepsilon,t,X^t \bigr) \bigr].
\]
We have
\begin{eqnarray*}
&&S (\vartheta,s,X_s ) -D(\hat\vartheta _\varepsilon,s,X_s)
\\
&&\quad=S(\vartheta,s,X_s)- S(\hat\vartheta_\varepsilon
,s,X_s)+S(\hat\vartheta_\varepsilon,s,X_s)
\\
&&\qquad{} - S\bigl(\hat\vartheta_\varepsilon ,s,x_s(\hat
\vartheta_\varepsilon)\bigr)- S'(\hat\vartheta_\varepsilon
,s,X_s) \bigl[X_s-x_s(\hat
\vartheta_\varepsilon) \bigr]
\\
&&\quad =- (\hat\vartheta_\varepsilon-\vartheta ) \dot S (\tilde
\vartheta,s,X_s )
\\
&&\qquad{} + \bigl[S'(\hat\vartheta_\varepsilon ,s,\tilde
X_s)-S'(\hat\vartheta_\varepsilon ,s,X_s)
\bigr] \bigl[X_s-x_s(\hat\vartheta_\varepsilon) \bigr]
\\
&&\quad =- (\hat\vartheta_\varepsilon-\vartheta ) \dot S (\tilde
\vartheta,s,X_s ) +\mathrm{O} \bigl(\varepsilon^2 \bigr).
\end{eqnarray*}
Therefore
\[
\varepsilon^{-1}Y_\varepsilon\bigl(\hat\vartheta_\varepsilon,t,X^t
\bigr) =-\frac{(\hat\vartheta_\varepsilon-\vartheta) }{\varepsilon
}\int_{0}^{t}
\frac{\dot S(\hat\vartheta_\varepsilon,s,X_s)^2
}{\sigma (s,X_s )^2} \,\mathrm{d}s+\mathrm{o} (1 ).
\]
Further,
\[
\varepsilon^{-1}\bigl(Z \bigl(\hat\vartheta_\varepsilon,t,X^t
\bigr)-Z \bigl(\vartheta,t,X^t \bigr) \bigr) =\frac{(\hat\vartheta_\varepsilon-\vartheta)}{\varepsilon} \dot
Z \bigl(\vartheta,t,X^t \bigr)+\mathrm{o} (1 ),
\]
where
\begin{eqnarray*}
\dot Z \bigl(\vartheta,t,X^t \bigr) &=& \int_{x_0}^{X_t}
\frac{\ddot S(\vartheta,t,y)}{ \sigma
 (t,y )^2} \,\mathrm{d}y - \int_{0}^{t}
\frac{\ddot
S (\vartheta,s,X_s )S (\vartheta
,s,X_s )}{ \sigma (s,X_s )^2} \,\mathrm{d}s
\\
&&{} -\int_{0}^{t}\int_{x_0}^{X_s}
\frac{ \ddot S_s'(\vartheta,s,y)
\sigma
 (s,y )-2\ddot S(\vartheta,s,y) \sigma_s'
 (s,y )}{ \sigma
 (s,y )^2} \,\mathrm{d}y \,\mathrm{d}s.
\end{eqnarray*}
We have uniform convergence of $X_t$ to $x_t$ w.r.t. $t$. Hence,
\[
\sup_{0\leq t\leq T}\bigl\llvert \dot Z \bigl(\vartheta,t,X^t
\bigr)-\dot Z \bigl(\vartheta,t,x^t \bigr)\bigr\rrvert \rightarrow0.
\]
Note that for any continuously differentiable function
$g (s,x )$ w.r.t. $s$ we have the relation
\[
\int_{x_0}^{x_t}g (t,y ) \,\mathrm{d}y - \int
_{0}^{t}g (s,x_s )S (\vartheta
,s,x_s ) \,\mathrm{d}s -\int_{0}^{t}
\int_{x_0}^{x_s} g_s'(s,y) \,
\mathrm{d}y \,\mathrm{d}s=0
\]
since
\[
\int_{0}^{t}g (s,x_s )S (\vartheta
,s,x_s ) \,\mathrm{d}s=\int_{0}^{t}g
(s,x_s )\,\mathrm{d}x_s
\]
and
\begin{eqnarray*}
\int_{0}^{t}g (t,x_s ) \,
\mathrm{d}x_s -\int_{0}^{t}g
(s,x_s ) \,\mathrm{d}x_s&=&\int_{0}^{t}
\int_{s}^{t} \frac{\partial g  (v,x_s )}{\partial v}\,\mathrm{d}v \,
\mathrm{d}x_s
\\
& =&\int_{0}^{t} \int_{0}^{t}
\one _{ \{v: x_v>x_s  \}} \frac{\partial g  (v,x_s )}{\partial v}\,\mathrm{d}v \,\mathrm{d}x_s
\\
&=&\int_{0}^{t}\int_{x_0}^{x_v}
g_v'(v,y) \,\mathrm{d}y \,\mathrm{d}v.
\end{eqnarray*}
Hence, $\dot Z (\vartheta,t,x^t )\equiv0$ for all
$t\in [0,T ]$.

By the It\^o formula,
\begin{eqnarray*}
\frac{Z (\vartheta,t,X^t )}{\varepsilon}&=&\frac{R
(\vartheta
,t,X^t )}{\varepsilon} -\int_{0}^{t}
\frac{\dot
S (\vartheta,s,X_s )S (\vartheta
,s,X_s )}{\varepsilon\sigma (s,X_s )^2} \,\mathrm{d}s
\\
&=&\int_{0}^{t}\frac{\dot
S (\vartheta,s,X_s )}{\varepsilon\sigma
(s,X_s )^2}
\,\mathrm{d}X_s-\int_{0}^{t}
\frac{\dot
S (\vartheta,s,X_s )S (\vartheta
,s,X_s )}{\varepsilon\sigma (s,X_s )^2} \,\mathrm{d}s
\\
&&{} +\frac{\varepsilon}{2}\int_{0}^{t}\sigma
(s,X_s )^2M_{xx}'' (
\vartheta,s,X_s )\,\mathrm{d}s
\\
&=&\int_{0}^{t}\frac{\dot
S (\vartheta,s,X_s )}{ \sigma (s,X_s )} \,
\mathrm{d}W_s+\mathrm{O} (\varepsilon ).
\end{eqnarray*}
Therefore, we obtain the convergence
\[
K_\varepsilon (t )\longrightarrow K (\vartheta ,t ).
\]
This convergence can be shown to be uniform w.r.t. $t$. This proves
the convergence $\delta_\varepsilon\rightarrow\delta$. Therefore the
Theorem~\ref{T1} is proved.
\end{pf}

Let us study the behaviour of the power function under the alternative. Suppose
that the observed diffusion process \eqref{0} has the trend coefficient
$S (t,x )$ which does not belong to the parametric
family. This family we described as follows:
\[
\mathcal{F}= \bigl\{S (\cdot )\dvtx S \bigl(\vartheta ,t,x_t (
\vartheta ) \bigr), 0\leq t\leq T, \vartheta\in \Theta \bigr\}.
\]
Here $x_t (\vartheta )$, $0\leq t\leq T$ is the solution of
equation \eqref{3}.

We introduce a slightly more strong condition of separability of the basic
hypothesis and the
alternative. Suppose that the function $S (t,x )$ satisfies
conditions $\mathcal{C}1$, $\mathcal{C}2$ and denote by $y_t$, $0\leq
t\leq T$ the
solution of the ordinary differential
equation obtained for $ (\varepsilon=0 )$
\[
\frac{\mathrm{d}y_t}{\mathrm{d}t}=S (t,y_t ),\qquad y_0=x_0.
\]
Then
\begin{eqnarray*}
\varepsilon^{-1} \bigl(X_t-x_t(\hat
\vartheta_\varepsilon ) \bigr)&=&\varepsilon^{-1}
(X_t-y_t )+\varepsilon ^{-1}
\bigl(y_t-x_t( \hat\vartheta_\varepsilon) \bigr)
\\
&=&y_t^{ (1 )}+ \varepsilon^{-1}
\bigl(y_t-x_t(\vartheta_*) \bigr)- \varepsilon^{-1}
(\hat\vartheta_\varepsilon-\vartheta_* ) \dot x_t(\vartheta_*)+
\mathrm{o} (1 ),
\end{eqnarray*}
where $y_t^{ (1 )}$ is a solution of the equation
\[
\mathrm{d}y_t^{ (1 )}=S' (t,y_t
)y_t^{
(1 )}\,\mathrm{d}t+\sigma (t,y_t )\,
\mathrm{d}W_t,\qquad y_0^{ (1 )}=0
\]
and $\vartheta_*$ is defined by the relation
%
%
\begin{equation}
\label{th} \inf_{\vartheta\in\Theta} \int_{0}^{T}
\biggl(\frac{S
(\vartheta,t
,y_t )-S (t,y_t )}{\sigma (t,y_t )} \biggr)^2\,\mathrm{d}t =\int
_{0}^{T} \biggl(\frac{S (\vartheta_*,t
,y_t )-S (t,y_t )}{\sigma (t,y_t )}
\biggr)^2\,\mathrm{d}t.
\end{equation}
Suppose that this equation has a unique solution $\vartheta_*$. Note that
$\varepsilon^{-1} (\hat\vartheta_\varepsilon-\vartheta_*
 ) $ is
tight (see Kutoyants \cite{r8} for details). Moreover, we
also suppose that the
basic
hypothesis and the alternative are separated in the following sense:
\[
\inf_{\vartheta\in\Theta} \int_{0}^{T}
\biggl(\frac{S
(\vartheta,t
,y_t )-S (t,y_t )}{\sigma (t,y_t )} \biggr)^2\,\mathrm{d}t >0.
\]

First, formally, we write
\begin{eqnarray*}
&&\int_{0}^{t} h_\varepsilon(\hat
\vartheta_\varepsilon,s)\,\mathrm{d}U_\varepsilon (s )
\\
&&\quad=\int_{0}^{t} \frac{\dot S(\hat
\vartheta
_\varepsilon,s,X_s)}{\sqrt{ \mathrm{I}(\hat\vartheta_\varepsilon)}
\sigma
 (s,X_s )} \,
\mathrm{d}W_s -\int_{0}^{t}
\frac{ \dot S(\hat\vartheta
_\varepsilon,s,X_s)  [S (s,X_s ) - D(\hat\vartheta
_\varepsilon,s,X_s) ]}{\sqrt{\mathrm{I}(\hat\vartheta
_\varepsilon) } \varepsilon\sigma
 (s,X_s )^2} \,\mathrm{d}s
\\
&&\quad =\int_{0}^{t} \frac{\dot S(\vartheta
_*,s,y_s)}{\sqrt{ \mathrm{I}(\vartheta_* )} \sigma
 (s,y_s )}\,
\mathrm{d}W_s-\int_{0}^{t}
\frac{ \dot
S(\vartheta
_*, s,X_s)  [S (s,X_s ) - S(\vartheta_*
,s,X_s) ]}{\sqrt{\mathrm{I}(\vartheta_* ) } \varepsilon\sigma
 (s,X_s )^2}\,\mathrm{d}s.
\end{eqnarray*}
Further
\begin{eqnarray*}
&&S (s,X_s )- D(\hat\vartheta_\varepsilon,s,X_s)
\\
&&\quad = S (s,X_s )-S \bigl(\hat\vartheta_\varepsilon,s,x_s
(\vartheta ) \bigr)-S' (\hat\vartheta_\varepsilon
,s,X_s ) \bigl(X_s-x_s (\hat
\vartheta_\varepsilon ) \bigr)
\\
&&\quad =S (s,X_s )-S (\hat\vartheta_\varepsilon
,s,X_s )+\mathrm{O} \bigl(\varepsilon^2 \bigr)
\\
&&\quad =S (s,X_s )-S (\vartheta_* ,s,X_s )+S (
\vartheta_*, s,X_s )- S (\hat \vartheta _\varepsilon,s,X_s
)+ \mathrm{O} \bigl(\varepsilon^2 \bigr)
\\
&&\quad =S (s,X_s )-S (\vartheta_* ,s,X_s )+ (\hat
\vartheta_\varepsilon- \vartheta_* )\dot S (\vartheta_*, s,X_s ) +
\mathrm{O} \bigl(\varepsilon^2 \bigr).
\end{eqnarray*}
Therefore,
\begin{eqnarray*}
\int_{0}^{t} h_\varepsilon(\hat
\vartheta_\varepsilon,s)\,\mathrm{d}U_\varepsilon (s ) &=&\int
_{0}^{t} \frac{\dot S(\vartheta
_*,s,y_s)}{\sqrt{ \mathrm{I}(\vartheta_* )} \sigma
 (s,y_s )}\,\mathrm{d}W_s-
\int_{0}^{t} \frac{(\hat\vartheta_\varepsilon- \vartheta_*)\dot S(\vartheta
_*,s,y_s)^2}{\varepsilon\sqrt{ \mathrm{I}(\vartheta_* )} \sigma
 (s,y_s )^2}\,\mathrm{d}s
\\
&&{} -\varepsilon^{-1}\int_{0}^{t}
\frac{ \dot S(\vartheta
_*, s,y_s)  [S (s,y_s ) - S(\vartheta_*
,s,y_s) ]}{\sqrt{\mathrm{I}(\vartheta_* ) } \sigma
 (s,y_s )^2}\,\mathrm{d}s+ \mathrm{O} \bigl(\varepsilon ^2 \bigr)
\\
& =&I_1 (t )-I_2 (t )-\varepsilon^{-1}I_3
(t )+ \mathrm{O} \bigl(\varepsilon ^2 \bigr)
\end{eqnarray*}
with an obvious notation. For the statistic $\delta_\varepsilon$ we
have the relations
%
%
\begin{equation}
\label{i3} \sqrt{\delta_\varepsilon}\geq\varepsilon^{-1}\bigl
\llVert I_3 (\cdot )h (\cdot )\bigr\rrVert -\bigl\llVert
I_1 (\cdot )h (\cdot )\bigr\rrVert -\bigl\llVert I_2 (
\cdot )h (\cdot )\bigr\rrVert +\mathrm{O} (\varepsilon ),
\end{equation}
where $h (\cdot )=h (\vartheta_*,s )$ and
$\llVert \cdot
\rrVert $ is the $L_2 (0,T )$ norm. Recall that the
quantities $\llVert  I_1 (\cdot
 )h (\cdot )\rrVert  $ and $\llVert  I_2 (\cdot
 )h (\cdot )\rrVert  $ are bounded in probability.

Introduce the condition

$\mathcal{C}5$. \textit{The functions $S (\vartheta,t,x
),S (t,x )$
and $\sigma (t,x )$ are such that}
\begin{eqnarray*}
&&\bigl\llVert I_3 (\cdot )h (\cdot )\bigr\rrVert ^2
\\
&&\quad =\int_{0}^{T} \biggl(\int
_{0}^{t} \frac{ \dot S(\vartheta
_*, s,y_s)  [S (s,y_s ) - S(\vartheta_*
,s,y_s) ]}{{\mathrm{I}(\vartheta_* ) } \sigma
 (s,y_s )^2}\,\mathrm{d}s
\biggr)^2\frac{\dot S
(\vartheta
_*,t )^2}{\sigma (t,y_t )^2}\,\mathrm{d}t >0. 
\end{eqnarray*}

This condition provides consistency of the test.

%
\begin{theorem}
\label{T2} Let conditions $\mathcal{C}1$--\hspace*{0.5pt}$\mathcal{C}5$ hold. Then the
test $\hat\psi
_\varepsilon$ is consistent.
\end{theorem}

\begin{pf}
The proof follows from the convergence $\delta
_\varepsilon
\rightarrow\infty$ under alternative (see \eqref{i3}).
\end{pf}

Note that if $\vartheta_*$ is an interior point of $\Theta$, then
\[
\int_{0}^{T} \frac{ \dot S(\vartheta
_*, s,y_s)  [S (s,y_s ) - S(\vartheta_*
,s,y_s) ]}{ \sigma
 (s,y_s )^2} \,\mathrm{d}s=0.
\]

If condition $\mathcal{C}5$ does not hold, then
\[
\int_{0}^{t} \frac{ \dot S(\vartheta
_*, s,y_s)  [S (s,y_s ) - S(\vartheta_*
,s,y_s) ]}{ \sigma
 (s,y_s )^2} \,\mathrm{d}s
\equiv0, \qquad\mbox{for all } t\in [0,T ].
\]
This equality is possible if
\[
\dot S(\vartheta_*, s,y_s) \bigl[S (s,y_s ) - S(
\vartheta_* ,s,y_s) \bigr]\equiv0, \qquad\mbox{for all } s\in [0,T ].
\]

An example of such \textit{invisible} alternative can be constructed as
follows: Suppose that the function $S (\vartheta,s,x )$
does not
depend on $\vartheta$ for $s\in [0,T/2 ]$, that is,
$S (\vartheta,s,x )=S_* (s,x )$ for all\vspace*{2pt}
$\vartheta\in
\Theta$. Then $\dot S(\vartheta_*, s,y_s)\equiv0 $ for $s\in
 [0,T/2 ]$. Therefore if $S (s,y_s ) =
S(\vartheta_*, s,y_s) $
for $s\in [T/2,T ]$ and a corresponding $\vartheta_*$ then
condition $\mathcal{C}5$ does not hold, but we can
have $S (s,y_s ) \neq S_*(s,y_s) $ for $s\in
[0,T/2 ]$.
This implies that the test $\hat\psi_\varepsilon$
is not consistent for this alternative.

\subsection{Second test}\label{sec3.2}

The second test is based on the following well-known transformation.
Suppose that
we have a Gaussian process $U(t)$, $0\leq t\leq1$ and $d\times d$
matrix $\NN (t ) $ defined by the relations
%
%
\begin{eqnarray}
\label{18} U (t )&=&W_t- \biggl\langle\int_{0}^{1}
\mathbf{h} (s )\,\mathrm{d}W_s,\int_{0}^{t}
\mathbf{h} (s ) \,\mathrm{d}s \biggr\rangle,
\\
\label{19}\NN (t )&=&\int_{t}^{1}\mathbf{h}
(s )\mathbf{h} (s )^*\,\mathrm{d} s,\qquad \NN (0 )=\JJ,
\end{eqnarray}
where $\JJ$ is the $d\times d$ unit matrix and $\mathbf{h}
(t )$ is a
continuous vector-valued function.

%
\begin{lemma}
\label{L2} Suppose that the matrix $\NN (t )$ is non-degenerate for all
$t\in[0,1)$. Then
%
%
\begin{equation}
\label{KK} U (t )+ \int_{0}^{t}\mathbf{h} (s
)^*\NN (s )^{-1}\int_{0}^{s}\mathbf{h}
(v ) \,\mathrm{d}U (v ) \,\mathrm{d}s=w (t ), \qquad 0\leq t\leq1,
\end{equation}
where $w (\cdot )$ is a Wiener process.
\end{lemma}

\begin{pf}
This formula was obtained by Khmaladze
\cite{r6}. The proof
there is based on two results: a result of Hitsuda \cite{r4} and
another one of Shepp \cite{r16}. Observe that there
are many publications
dealing with this
transformation (see, e.g., the paper
Maglaperidze \textit{et al.}
\cite{r15} and the references
therein). Another direct proof is given in Kleptsyna and Kutoyants \cite{r7}.
\end{pf}

Note that representation \eqref{18} and \eqref{19} implies that
%
%
\begin{equation}
\label{21} \int_{0}^{1}\mathbf{h} (s )\,
\mathrm{d}U (s )=0.
\end{equation}

Suppose that $\vartheta\in\Theta
$. Here $\Theta$ is an open bounded subset of $\mathscr{R}^d$. Now
$\mathbf{h} (\vartheta,s )$, $\mathbf{R} (\vartheta
,t,X^t ) $ and
$\mathbf{Q} (\vartheta,t,X^t )$ are $d$-vectors and the
Fisher information
$\II (\vartheta ) $ is a $d\times d$ matrix.

Introduce the following stochastic processes:
\begin{eqnarray*}
\bar\mathbf{h}_\varepsilon(\vartheta,t )&=&\frac{\dot\mathbf{S}(\vartheta,t,X_t)}{\sigma (t,X_t
)},
\\
\bar\NN (\vartheta,t)&=&\int_{t}^{T}
\frac{\dot\mathbf{S}(\vartheta
,s,x_s)\dot
\mathbf{S}(\vartheta,s,x_s)^* }{\sigma(s,x_s)^2} \,\mathrm{d} s,
\\
\bar\NN_\varepsilon (\vartheta,t)&=&\int_{t}^{T}
\frac{\dot\mathbf{S}(\vartheta
,s,X_s)\dot
\mathbf{S}(\vartheta,s,X_s)^* }{\sigma(s,X_s)^2}\,\mathrm{d} s,
\end{eqnarray*}
and put
\[
\Delta_\varepsilon=\frac{1}{T^2}\int_{0}^{T}W_\varepsilon
(t )^2\,\mathrm{d}t.
\]
Here
%
%
\begin{eqnarray}
\label{We} W_\varepsilon (t )&=&\int_{0}^{t}
\frac{\mathrm{d}X_s}{\varepsilon\sigma
 (s,X_s )}- \int_{0}^{t}\frac{D(\hat\vartheta
_\varepsilon
,s,X_s)}{\varepsilon\sigma (s,X_s ) }
\,\mathrm{d}s
\nonumber
\\[-8pt]
\\[-8pt]
&&{}+\varepsilon^{-1}\int_{0}^{t}\bar
\mathbf{h}_\varepsilon(\hat \vartheta_\varepsilon ,s)^* \bar
\NN_\varepsilon (\hat\vartheta_\varepsilon ,s )_+^{-1} \bigl[
\mathbf{R} \bigl(\hat\vartheta_\varepsilon ,s,X^s \bigr) -
\mathbf{Q} \bigl(\hat\vartheta_\varepsilon,s,X^s \bigr) \bigr] \,
\mathrm{d}s.
\nonumber
\end{eqnarray}
We use the following convention for the matrix $\bar\NN$:
\[
\bar\NN^{-1}_+= %
\cases{ \bar\NN^{-1},&\quad$
\mbox{if } \bar\NN \mbox{ is non-degenerate}$,
\cr
0, &\quad$\mbox{if } \bar\NN
\mbox{ is degenerate}$. } %
\]
We have the following result.
%

\begin{theorem}\label{T3}
Suppose that conditions $\mathcal{C}2$--\hspace*{0.5pt}$\mathcal{C}4$ hold and
the matrix $\bar\NN
(\vartheta,t)$ is uniformly in $\vartheta\in\Theta$ non-degenerate for
all $t\in[0,1)$. Then the test
\[
\hat\Psi_\varepsilon=\one _{ \{\Delta_\varepsilon>c_\alpha
 \}}, \qquad\Pb \biggl(\int
_{0}^{1}w (s )^2\,
\mathrm{d}s>c_\alpha \biggr)=\alpha
\]
is ADF and belongs to $\mathcal{K}_\alpha$.
\end{theorem}

\begin{pf} We have to show that under hypothesis $\mathcal{H}_0$ the
convergence
%
%
\begin{equation}
\label{q} \Delta_\varepsilon\Longrightarrow\Delta=\int
_{0}^{1}w (s )^2\,\mathrm{d}s
\end{equation}
holds.

The construction of the ADF GoF test is based on Lemmas~\ref{L1} and~\ref{L2}. We have the similar to \eqref{8} presentation \eqref{18} with
$\mathbf{h} (\vartheta,t )$ defined in \eqref{h}. Let us
denote $U_\varepsilon (\cdot ), \mathbf{h}_\varepsilon
(\hat\vartheta_\varepsilon
,\cdot)$, and $\NN_\varepsilon (\cdot )$ the \textit
{empirical versions} of
$U (\cdot ),\mathbf{h} (\vartheta,\cdot
 )$ and
\[
\NN (\vartheta,t ) =\II(\vartheta)^{-1}\int_{t}^{T}
\frac{\dot\mathbf{S}(\vartheta
,s,x_s)\dot
\mathbf{S}(\vartheta,s,x_s)^* }{\sigma(s,x_s)^2}\, \mathrm{d} s,\qquad \NN (\vartheta,0 ) =\JJ,
\]
respectively:
\begin{eqnarray*}
U_\varepsilon (t )&=&\int_{0}^{t}
\frac{\psi_\varepsilon (s
) }{\sigma
 (s,X_s )} \,\mathrm{d} V_\varepsilon (s ),
\\
V_\varepsilon (t )&=&\frac{X_t-x_t(\hat\vartheta
_\varepsilon)}{\psi
_\varepsilon (t )\varepsilon},
\\
\mathbf{h}_\varepsilon(\hat\vartheta_\varepsilon,t )&=&\II_\varepsilon(
\hat\vartheta_\varepsilon)^{-1/2}\frac{\dot
\mathbf{S}(\hat\vartheta
_\varepsilon,t,X_t)}{\sigma (t,X_t )},
\\
\II_\varepsilon(\hat\vartheta_\varepsilon)&=&\int_{0}^{T}
\frac
{\dot\mathbf{S}(\hat\vartheta
_\varepsilon,s,X_s)\dot
\mathbf{S}(\hat\vartheta_\varepsilon,s,X_s)^* }{\sigma(s,X_s)^2} \,\mathrm{d} s,
\\
\NN_\varepsilon (\hat\vartheta_\varepsilon,t)&=&\II_\varepsilon(\hat
\vartheta _\varepsilon )^{-1}\int_{t}^{T}
\frac{\dot\mathbf{S}(\hat\vartheta_\varepsilon
,s,X_s)\dot
\mathbf{S}(\hat\vartheta_\varepsilon,s,X_s)^* }{\sigma(s,X_s)^2} \,\mathrm{d} s.
\end{eqnarray*}
Recall that there is a problem of definition of the integral for
$U_\varepsilon (\cdot )$ because the integrand depends on the
future. As convergence is uniform w.r.t. $t\in [0,T-\nu ]$:
\[
\mathbf{h}_\varepsilon(\hat\vartheta_\varepsilon,t )\longrightarrow
\mathbf{h} (\vartheta ,t ),\qquad \II_\varepsilon(\hat\vartheta_\varepsilon
)\longrightarrow \II(\vartheta),\qquad \NN_\varepsilon(\hat
\vartheta_\varepsilon ,t)\longrightarrow\NN(\vartheta,t).
\]
The required limits can be obtained.

Introduce (formally) the statistic
%
%
\begin{equation}
\label{W} W^\star_\varepsilon (t )=U_\varepsilon (t )+\int
_{0}^{t}\mathbf{h}_\varepsilon(\hat\vartheta
_\varepsilon ,s)^* \NN_\varepsilon (\hat\vartheta_\varepsilon,s
)_+^{-1} \int_{0}^{s}
\mathbf{h}_\varepsilon (\hat\vartheta_\varepsilon,v) \,\mathrm{d}U_\varepsilon
(v ) \,\mathrm{d}s.
\end{equation}

Observe that
\begin{eqnarray*}
&&\mathbf{h} (\vartheta,s )^* \NN (\vartheta,s )^{-1}\mathbf{h} (
\vartheta ,v )
\\
&&\quad=\frac{\dot\mathbf{S} (\vartheta,s,x_s )^*}{\sigma
 (s,x_s )} \biggl(\int_{s}^{T}
\frac{\dot\mathbf{S}
(\vartheta
,r,x_r )\dot\mathbf{S} (\vartheta,r,x_r )^*
}{\sigma
 (r,x_r )^2} \,\mathrm{d} r \biggr)^{-1} \frac{\dot\mathbf{S} (\vartheta
,v,x_v )}{\sigma (v,x_v )}.
\end{eqnarray*}
Therefore this term does not depend on the information matrix $ \II
 (\vartheta
 ) $ and we can replace the statistics $\mathbf{h}_\varepsilon
(\hat\vartheta_\varepsilon,s) $ and $\NN_\varepsilon
 (\hat\vartheta_\varepsilon,s ) $ in \eqref{W} by $\bar
\mathbf{h}_\varepsilon
(\hat\vartheta_\varepsilon,s) $ and $\bar\NN_\varepsilon
 (\hat\vartheta_\varepsilon,s ) $.

For the process $U_\varepsilon (\cdot )$, we have equality
\eqref{u} (formally)
\[
U_\varepsilon (t )=\int_{0}^{t}
\frac{\mathrm{d}X_s}{\varepsilon\sigma
 (s,X_s )}- \int_{0}^{t}\frac{D(\hat\vartheta
_\varepsilon
,s,X_s)}{\varepsilon\sigma
 (s,X_s ) }
\,\mathrm{d}s.
\]

Hence, we obtain the vector-valued integral
\[
\int_{0}^{t}\bar\mathbf{h}_\varepsilon(\hat
\vartheta_\varepsilon ,s)\,\mathrm{d}U_\varepsilon (s ) =\int
_{0}^{t} \frac{\dot\mathbf{S}(\hat\vartheta
_\varepsilon,s,X_s)}{\varepsilon\sigma
 (s,X_s )^2}\,
\mathrm{d}X_s-\int_{0}^{t}
\frac{\dot\mathbf{S}(\hat\vartheta
_\varepsilon
,s,X_s)D(\hat\vartheta_\varepsilon,s,X_s)}{\varepsilon\sigma
 (s,X_s )^2}\,\mathrm{d}s.
\]
Introduce the vector-function
\[
\mathbf{M} (\vartheta,t,x )=\int_{x_0}^{x}
\frac{\dot
\mathbf{S}(\vartheta,t,y)}{ \sigma
 (t,y )^2} \,\mathrm{d}y.
\]
Then by the It\^o formula
\begin{eqnarray*}
\int_{0}^{t} \frac{\dot\mathbf{S}(\vartheta,s,X_s)}{\sigma
 (s,X_s )^2}\,
\mathrm{d}X_s&=&\int_{x_0}^{X_t}
\frac{\dot
\mathbf{S}(\vartheta,t,y)}{ \sigma
 (t,y )^2} \,\mathrm{d}y-\int_{0}^{t}\int
_{x_0}^{X_s} \frac
{\dot\mathbf{S}_s'(\vartheta,s,y)}{ \sigma
 (s,y )^2}\, \mathrm{d}s
\\
&&{} + \int_{0}^{t}\int_{x_0}^{X_s}
\frac{2\dot\mathbf{S}(\vartheta
,s,y)\sigma_s'
 (s,y )}{ \sigma
 (s,y )^3} \,\mathrm{d}s+\mathrm{O} \bigl(\varepsilon^2
\bigr).
\end{eqnarray*}
Put\vspace*{1pt}
\[
\mathbf{K}_\varepsilon (t )=\int_{0}^{t}\bar
\mathbf{h}_\varepsilon(\hat\vartheta_\varepsilon,s)\,\mathrm{d}U_\varepsilon
(s ) =\varepsilon^{-1} \bigl[\mathbf{R} \bigl(\hat\vartheta_\varepsilon
,t,X^t \bigr)-\mathbf{Q} \bigl(\hat\vartheta_\varepsilon
,t,X^t \bigr) \bigr].
\]
Note that we have dropped the term of order $\mathrm{O} (\varepsilon
^2 )$.

Then formal expression \eqref{W} for
$W^\star_\varepsilon (t )$ can be replaced by \eqref{We}\vspace*{1pt}
\begin{eqnarray*}
W_\varepsilon (t )&=& \int_{0}^{t}
\frac{\mathrm{d}X_s}{\varepsilon\sigma
 (s,X_s )}- \int_{0}^{t}\frac{D(\hat\vartheta
_\varepsilon
,s,X_s)}{\varepsilon\sigma (s,X_s ) }
\,\mathrm{d}s
\\
&&{}+\varepsilon^{-1}\int_{0}^{t}\bar
\mathbf{h}_\varepsilon(\hat \vartheta_\varepsilon ,s)^* \bar
\NN_\varepsilon (\hat\vartheta_\varepsilon ,s )_+^{-1} \bigl[ R
\bigl(\hat\vartheta_\varepsilon ,s,X^s \bigr) -Q \bigl(\hat
\vartheta_\varepsilon,s,X^s \bigr) \bigr] \,\mathrm{d}s.
\end{eqnarray*}
For the first two terms of $W_\varepsilon (t )$
we have
\begin{eqnarray*}
U_\varepsilon (t )&=&\int_{0}^{t}
\frac{\mathrm{d}X_s}{\varepsilon\sigma
 (s,X_s )}- \int_{0}^{t}\frac{D(\hat\vartheta
_\varepsilon
,s,X_s)}{\varepsilon\sigma (s,X_s ) }
\,\mathrm{d}s
\\
& =&W_t+\int_{0}^{t}
\frac{S (\vartheta,s,X_s ) -
S(\hat\vartheta_\varepsilon
,s,x_s(\hat\vartheta_\varepsilon))-S'(\hat\vartheta_\varepsilon
,s,X_s)(X_s-x_s(\hat\vartheta
_\varepsilon))}{\varepsilon\sigma (s,X_s ) }\,\mathrm{d}s
\\
& =&W_t- \biggl\langle\frac{\hat\vartheta_\varepsilon-\vartheta
}{\varepsilon
},\int_{0}^{t}
\frac{\dot S(\tilde\vartheta,s,X_s)}{\sigma
 (s,X_s ) }\,\mathrm{d}s \biggr\rangle
\\
&&{} +\int_{0}^{t}\frac{ [S'(\hat\vartheta_\varepsilon
,s,\tilde X_s)-S'(\hat\vartheta_\varepsilon
,s,X_s) ](X_s-x_s(\hat\vartheta
_\varepsilon))}{\varepsilon\sigma (s,X_s ) }\,
\mathrm{d}s
\\
& =&W_t- \biggl\langle\mathrm{I} (\vartheta )^{-1}\int
_{0}^{T}\frac{\dot\mathbf{S}( \vartheta,s,x_s)}{
\sigma
 (s,x_s ) }\,\mathrm{d}W_s,
\int_{0}^{t}\frac{\dot\mathbf{S}( \vartheta
,s,x_s)}{ \sigma
 (s,x_s ) }\,\mathrm{d}s
\biggr\rangle+\mathrm{o} (1 )
\\
& =&U (\vartheta,t )+\mathrm{o} (1 ).
\end{eqnarray*}
Here $\llvert \tilde\vartheta-\vartheta\rrvert \leq\llvert \hat
\vartheta
_\varepsilon\rrvert $ and
\begin{eqnarray*}
\llvert \tilde X_s-X_s\rrvert &\leq& \bigl\llvert
x_s(\hat\vartheta_\varepsilon)-X_s\bigr\rrvert
\\[1pt]
&\leq&
\bigl\llvert x_s(\hat\vartheta_\varepsilon)-x_s (
\vartheta )\bigr\rrvert +\bigl\llvert x_s (\vartheta )-X_s
\bigr\rrvert \rightarrow0.
\end{eqnarray*}
This convergence is uniform w.r.t. $s\in [0,T ]$. Hence,
\[
\sup_{0\leq t\leq T}\bigl\llvert U_\varepsilon (t )-U (\vartheta,t )
\bigr\rrvert \longrightarrow0.
\]
Further, similar arguments give the uniform convergence
w.r.t. $t\in [0,T ]$
\[
\bar\mathbf{h}_\varepsilon(\hat\vartheta_\varepsilon,t)=
\frac
{\dot
\mathbf{S}(\hat\vartheta_\varepsilon,t,X_t)}{\sigma
(t,X_t )} \rightarrow\bar\mathbf{h} (\vartheta,t ), \qquad \bar\NN
_\varepsilon (\hat\vartheta_\varepsilon,t)\rightarrow\bar\NN (\vartheta,t).
\]
We have to show that $\mathbf{K}_\varepsilon (t )
\longrightarrow
\mathbf{K} (\vartheta,t ) $, where
\[
\mathbf{K} (\vartheta,t )= \int_{0}^{t} \bar
\mathbf{h} (\vartheta ,s )\,\mathrm{d}W_s- \int_{0}^{T}
\bar\mathbf{h} (\vartheta ,s )\,\mathrm{d}W_s\int
_{0}^{t} \bar\mathbf{h} (\vartheta ,s )\bar
\mathbf{h} (\vartheta,s )^*\,\mathrm{d}s.
\]
Denote
\begin{eqnarray*}
\mathbf{Y}_\varepsilon \bigl(\hat\vartheta_\varepsilon,t,X^t
\bigr)&=&\int_{0}^{t}\frac{\dot\mathbf{S} (\hat\vartheta
_\varepsilon
,s,X_s ) [ S (\vartheta,s,X_s ) -D(\hat
\vartheta
_\varepsilon,s,X_s)  ]}{\sigma (s,X_s )^2} \,
\mathrm{d}s,
\\
\mathbf{Z} \bigl(\hat\vartheta_\varepsilon,t,X^t \bigr)&=&{
\mathbf{R} \bigl(\hat\vartheta _\varepsilon,t,X^t \bigr)}- \int
_{0}^{t}\frac{\dot\mathbf{S}
(\hat\vartheta
_\varepsilon,s,X_s )S (\vartheta,s,X_s )}{\sigma
 (s,X_s )^2} \,\mathrm{d}s.
\end{eqnarray*}
Then
\[
\mathbf{K}_\varepsilon (t )=\varepsilon^{-1} \bigl[
\mathbf{Y}_\varepsilon \bigl(\hat\vartheta_\varepsilon,t,X^t
\bigr)+ \mathbf{Z} \bigl(\hat\vartheta_\varepsilon,t,X^t \bigr)
\bigr].
\]
We have
\begin{eqnarray*}
&& S (\vartheta,s,X_s ) -D(\hat\vartheta _\varepsilon,s,X_s)
 \\
 &&\quad= S(\vartheta,s,X_s)- S(\hat\vartheta_\varepsilon
,s,X_s)+S(\hat\vartheta_\varepsilon,s,X_s)
\\
&&\qquad {} - S\bigl(\hat\vartheta_\varepsilon ,s,x_s(\hat
\vartheta_\varepsilon)\bigr)- S'(\hat\vartheta_\varepsilon
,s,X_s) \bigl[X_s-x_s(\hat
\vartheta_\varepsilon) \bigr]
\\
&&\quad  =- \bigl\langle (\hat\vartheta_\varepsilon-\vartheta ), \dot
\mathbf{S} (\tilde\vartheta,s,X_s ) \bigr\rangle
\\
&&\qquad {} + \bigl[S'(\hat\vartheta_\varepsilon ,s,\tilde
X_s)-S'(\hat\vartheta_\varepsilon ,s,X_s)
\bigr] \bigl[X_s-x_s(\hat\vartheta_\varepsilon) \bigr]
\\
&&\quad  =- \bigl\langle (\hat\vartheta_\varepsilon-\vartheta ), \dot
\mathbf{S} (\tilde\vartheta,s,X_s ) \bigr\rangle+\mathrm{O} \bigl(
\varepsilon^2 \bigr).
\end{eqnarray*}
Therefore
\[
\varepsilon^{-1}\mathbf{Y}_\varepsilon \bigl(\hat\vartheta
_\varepsilon,t,X^t \bigr) =-\frac{ (\hat\vartheta_\varepsilon-\vartheta )
}{\varepsilon
}\int
_{0}^{t}\frac{\dot\mathbf{S} (\hat\vartheta_\varepsilon
,s,X_s )\dot
\mathbf{S} (\hat\vartheta_\varepsilon,s,X_s )^* }{\sigma
 (s,X_s )^2} \,\mathrm{d}s.
\]

Further,
\[
\varepsilon^{-1} \bigl({\mathbf{Z} \bigl(\hat\vartheta
_\varepsilon,t,X^t \bigr)-\mathbf{Z} \bigl(\vartheta,t,X^t
\bigr)} \bigr) =\frac{\hat\vartheta_\varepsilon-\vartheta
}{\varepsilon} \dot \ZZ \bigl(\vartheta,t,X^t
\bigr)+\mathrm{o} (1 ),
\]
where
\begin{eqnarray*}
\dot\ZZ \bigl(\vartheta,t,X^t \bigr) &=& \int_{x_0}^{X_t}
\frac{\ddot\SS(\vartheta,t,y)}{ \sigma
 (t,y )^2} \,\mathrm{d}y - \int_{0}^{t}
\frac{\ddot
\SS (\vartheta,s,X_s )S (\vartheta
,s,X_s )}{ \sigma (s,X_s )^2} \,\mathrm{d}s
\\
&&{} -\int_{0}^{t}\int_{x_0}^{X_s}
\frac{ \ddot\SS_s'(\vartheta
,s,y)\sigma
 (s,y )-2\ddot\SS(\vartheta,s,y)\sigma_s'
 (s,y )}{ \sigma
 (s,y )^2} \,\mathrm{d}y \,\mathrm{d}s.
\end{eqnarray*}
Here $\ddot\SS (\cdot )$ is the matrix of second derivatives
w.r.t. $\vartheta$. We have uniform convergence of $X_t$ to $x_t$
w.r.t. $t$, hence
\[
\sup_{0\leq t\leq T}\bigl\llvert \dot\ZZ \bigl(\vartheta,t,X^t
\bigr)-\dot\ZZ \bigl(\vartheta,t,x^t \bigr)\bigr\rrvert \rightarrow0.
\]
Observe that for any continuously differentiable function
$g (s,x )$ w.r.t. $s$ we have
\[
\int_{x_0}^{x_t}g (t,y ) \,\mathrm{d}y - \int
_{0}^{t}g (s,x_s )S (\vartheta
,s,x_s ) \,\mathrm{d}s -\int_{0}^{t}
\int_{x_0}^{x_s} g_s'(s,y) \,
\mathrm{d}y \,\mathrm{d}s=0
\]
since
\[
\int_{0}^{t}g (s,x_s )S (\vartheta
,s,x_s ) \,\mathrm{d}s=\int_{0}^{t}g
(s,x_s )\,\mathrm{d}x_s
\]
and
\begin{eqnarray*}
 && \int_{0}^{t}g (t,x_s )\,
\mathrm{d}x_s -\int_{0}^{t}g
(s,x_s ) \,\mathrm{d}x_s
 \\
&& \quad=\int_{0}^{t} \int_{s}^{t}
\frac{\partial g  (v,x_s )}{\partial v}\,\mathrm{d}v \,\mathrm{d}x_s
\\
&&\quad  =\int_{0}^{t} \int_{0}^{t}
\one _{ \{v: x_v>x_s  \}} \frac{\partial g  (v,x_s )}{\partial v}
\,\mathrm{d}v\, \mathrm{d}x_s =\int
_{0}^{t}\int_{x_0}^{x_v}
g_v'(v,y) \,\mathrm{d}y \,\mathrm{d}v.
\end{eqnarray*}
Hence, $\dot\ZZ (\vartheta,t,x^t )\equiv0$ for all
$t\in [0,T ]$.

By the It\^o formula
\begin{eqnarray*}
\frac{\mathbf{Z} (\vartheta,t,X^t )}{\varepsilon}&=&\frac
{\mathbf{R} (\vartheta
,t,X^t )}{\varepsilon} -\int_{0}^{t}
\frac{\dot
\mathbf{S} (\vartheta,s,X_s )S (\vartheta
,s,X_s )}{\varepsilon\sigma (s,X_s )^2} \,\mathrm{d}s
\\
&=&\int_{0}^{t}\frac{\dot
\mathbf{S} (\vartheta,s,X_s )}{\varepsilon\sigma
(s,X_s )^2}\,
\mathrm{d}X_s-\int_{0}^{t}
\frac{\dot
\mathbf{S} (\vartheta,s,X_s )S (\vartheta
,s,X_s )}{\varepsilon\sigma (s,X_s )^2} \,\mathrm{d}s
\\
&&{} +\frac{\varepsilon}{2}\int_{0}^{t}\sigma
(s,X_s )^2\mathbf{M}_{xx}''
(\vartheta,s,X_s )\,\mathrm{d}s
\\
&=&\int_{0}^{t}\frac{\dot
\mathbf{S} (\vartheta,s,X_s )}{ \sigma (s,X_s
)} \,
\mathrm{d}W_s+\mathrm{O} (\varepsilon ).
\end{eqnarray*}
Therefore, we obtain the convergence
\begin{eqnarray*}
\mathbf{K}_\varepsilon (t )&=&\varepsilon^{-1} \bigl(\mathbf{R}
\bigl(\hat\vartheta_\varepsilon ,t,X^t \bigr)- \mathbf{Q} \bigl(
\hat\vartheta_\varepsilon ,t,X^t \bigr) \bigr)
\\
&=&\varepsilon^{-1} \bigl(\mathbf{Y} \bigl(\hat\vartheta_\varepsilon
,t,X^t \bigr)+\mathbf{Z} \bigl(\hat\vartheta_\varepsilon
,t,X^t \bigr) \bigr)\longrightarrow\mathbf{K} (\vartheta ,t ).
\end{eqnarray*}
Further, the matrix $\bar\NN_\varepsilon(\hat\vartheta_\varepsilon
,s)$ converges uniformly in $s\in [0,T ]$ to the matrix
$\bar
N (\vartheta,s )$. Therefore, for $\nu>0$ we have uniform
on $s\in
 [0,T-\nu ]$ convergence of $\bar\NN_\varepsilon
(\hat\vartheta_\varepsilon,s)_+^{-1}$ to $N (\vartheta,s )^{-1}
$. Introduce the random function
\[
y_\varepsilon (s )=\varepsilon^{-1}\bar\mathbf{h}_\varepsilon(
\hat\vartheta_\varepsilon ,s)^* \bar\NN_\varepsilon (\hat
\vartheta_\varepsilon ,s )_+^{-1} \bigl[ \mathbf{R} \bigl(\hat
\vartheta_\varepsilon ,s,X^s \bigr) -\mathbf{Q} \bigl(\hat
\vartheta_\varepsilon,s,X^s \bigr) \bigr].
\]
It is shown that we have convergence
\[
\sup_{0\leq s\leq T-\nu}\bigl\llvert y_\varepsilon (s )-y (\vartheta,s )
\bigr\rrvert \longrightarrow0,
\]
where
\[
y (\vartheta,s )=\bar\mathbf{h} (\vartheta,s)^* \bar\NN (\vartheta ,s
)^{-1}\mathbf{K} (\vartheta,s ).
\]

Hence we also have convergence for all $t\in[0,1)$
\[
W_\varepsilon (t )\longrightarrow U (\vartheta ,t )+\int_{0}^{t}
\bar\mathbf{h} (\vartheta,s )^*\bar N (\vartheta ,s )^{-1}\mathbf{K} (
\vartheta,s )\,\mathrm{d}s=w (t ).
\]
A similar argument can show that for any $0\leq t_1<\cdots<t_k\leq T$
we have convergence of the vectors
\[
\bigl(W_\varepsilon (t_1 ),\ldots,W_\varepsilon
(t_k ) \bigr)\Longrightarrow \bigl( w (t_1 ),\ldots,w
(t_k ) \bigr).
\]
Further, a direct but cumbersome calculation allows us to write the estimate
\[
\Ex_\vartheta\bigl\llvert W_\varepsilon (t_1
)-W_\varepsilon (t_2 ) \bigr\rrvert ^2\leq C\llvert
t_2-t_1\rrvert , \qquad t_1,t_2
\in [0,T-\nu ].
\]
These two conditions provide weak convergence of the integrals
\[
\int_{0}^{T-\nu}W_\varepsilon (t )^2
\,\mathrm{d}t\Longrightarrow \int_{0}^{T-\nu}w (t
)^2\,\mathrm{d}t
\]
for any $\nu>0$. It can be shown that for any $\eta>0$ there exist
$\nu>0$
such that
\[
\int_{T-\nu}^{T}\Ex_\vartheta W_\varepsilon
(t )^2\,\mathrm{d}t\leq\eta.
\]
The proof is close to that given in Maglaperidze
\textit{et al.}
\cite{r15}
for a similar integral.
\end{pf}

\section{Examples}

\begin{Example}\label{exa1}
We consider the simplest case which allows us to
have an ADF
GoF test for each~$\varepsilon$, that is, no need to study statistics
as $\varepsilon
\rightarrow0$. Observe that a similar
situation is discussed in Khmaladze \cite{r6} but for a
different problem.

Suppose that the observed diffusion process (under hypothesis)
is
%
%
\begin{equation}
\label{e1} \mathrm{d}X_t=\vartheta\,\mathrm{d}t+\varepsilon\,
\mathrm{d}W_t,\qquad X_0=0, 0\leq t\leq1.
\end{equation}
Then
\begin{eqnarray*}
h (\vartheta,t )&=&1,\qquad \mathrm{I} (\vartheta )=1,\qquad N (\vartheta,t )
=1-t,
\\
\hat\vartheta_\varepsilon& =&X_1,\qquad \varepsilon^{-1}
(\hat \vartheta _\varepsilon-\vartheta ) =W_1\sim\mathcal{N} (0,1
).
\end{eqnarray*}
Further
\begin{eqnarray*}
x_t (\vartheta )&=&\vartheta t,\qquad x_t^{ (1 )}
(\vartheta )= W_t,\qquad U (\vartheta ,t )=W_t-W_1
t,
\\
V_\varepsilon (t )&=&U_\varepsilon (t )=\varepsilon^{-1}
(X_t-X_1 t )=W_t-W_1 t=B (t ).
\end{eqnarray*}
Therefore,
\[
W_\varepsilon (t )=\varepsilon^{-1} (X_t-X_1
t )+\varepsilon^{-1}\int_{0}^{t} (1-s
)^{-1} [X_s-X_1 s ] \,\mathrm{d}s
\]
and under the basic hypothesis we have
\[
W_\varepsilon (t )= B (t )+\int_{0}^{t}
\frac{B (s )}{1-s} \,\mathrm{d}s =w (t ).
\]
Therefore,
\[
\Delta_\varepsilon=\int_{0}^{1}W_\varepsilon
(t )^2\,\mathrm{d}t=\int_{0}^{1}w (t
)^2\,\mathrm{d}t
\]
and the test
$\hat\Psi_\varepsilon=\one _{ \{\Delta_\varepsilon>c_\alpha
 \}}\in
\mathcal{K}_\alpha$
satisfies the equality
\[
\Ex_\vartheta\hat\Psi_\varepsilon=\Pb \biggl\{\int_{0}^{1}w
(t )^2\,\mathrm{d}t>c_\alpha \biggr\}=\alpha.
\]
\end{Example}

\begin{Example}\label{exa2}
Consider the linear case
\[
\mathrm{d}X_t=\bigl\langle\vartheta,\mathbf{H} (t,X_t )
\bigr\rangle\, \mathrm{d}t+\varepsilon \sigma (t,X_t ) \,
\mathrm{d}W_t,\qquad X_0=x_0, 0\leq t\leq T,
\]
where $\vartheta\in\Theta\subset\mathcal{R}^d$ and assume that the
functions
$\mathbf{H} (t,x )$ and $\sigma (t,x )$ satisfy
regularity
conditions. In this case, we can take $\bar\mathbf{h}_\varepsilon
 (\vartheta
,t )=\bar\mathbf{h}_\varepsilon (t ) $, that is,
this vector-valued function
does not depend on $\vartheta$. Hence, the stochastic integral is well defined
and the test has a simplified form. We have
\begin{eqnarray*}
\bar\mathbf{h}_\varepsilon (t )&=&\frac{\mathbf{
H} (t,X_t )}{\sigma (t,X_t )},\qquad \bar{
\NN}_\varepsilon (\vartheta,s )=\int_{s}^{T}
\frac
{\mathbf{
H} (t,x_t (\vartheta ) )\mathbf{
H} (t,x_t (\vartheta ) )^*}{\sigma
 (t,x_t (\vartheta ) )^2} \,\mathrm{d}s,
\\
\mathrm{d}U_\varepsilon (t )&=&\frac{\mathrm{d}X_t}{\varepsilon\sigma\varphi
 (t,X_t )}- \frac{ [\langle\hat\vartheta
_\varepsilon,
\mathbf{H}(t, x_t(\hat\vartheta_\varepsilon))\rangle+\langle\hat
\vartheta
_\varepsilon, \mathbf{H}_x'(t,X_t)\rangle(X_t-x_t(\hat\vartheta
_\varepsilon
)) ]\,\mathrm{d}t }{\varepsilon\sigma (t,X_t )},
\\
W_\varepsilon (t )&=&U_\varepsilon (t )+\int_{0}^{t}
\mathbf{H} (s,X_s )^* \bar{\NN}_\varepsilon (\hat
\vartheta_\varepsilon,s ) ^{-1}\int_{0}^{s}
\mathbf{H} (v,X_v ) \,\mathrm{d}U_\varepsilon (v ) \,\mathrm{d}s
\end{eqnarray*}
and so on.
\end{Example}


\section*{Acknowledgements}

This study was partially supported by Russian Science Foundation (research
project No. 14-49-00079). The author thanks the referee for helpful comments.


%
%

\printhistory
\end{document}